\documentclass[11pt]{amsart}
\usepackage{amssymb}
\usepackage{verbatim}  

\newtheorem{thm}	  {Theorem}      [section]
\newtheorem{prop}   [thm] {Proposition}
\newtheorem{lemma}  [thm] {Lemma}
\newtheorem{cor}    [thm] {Corollary}

\theoremstyle{remark}
\newtheorem*{remark}{Remark}

\newcommand\F{\mathbb F}
\newcommand\Fq{\F_q}
\newcommand\Pone{\mathbb P^1}
\DeclareMathOperator\T{\mathbb T}
\DeclareMathOperator\Aut{Aut}
\DeclareMathOperator\Gal{Gal}
\DeclareMathOperator\PSL{PSL_2}
\DeclareMathOperator\PGL{PGL_2}
\DeclareMathOperator\PGammaL{P\Gamma L_2}
\DeclareMathOperator\SL{SL_2}

\DeclareMathOperator\ind{ind}
\DeclareMathOperator\Sym{Sym}
\renewcommand{\bar}[1]{#1\llap{$\overline{\phantom{\rm#1}}$}}
\newcommand{\col}{\,{:}\,}
\newcommand\bl{\bar \ell}
\newcommand\bbl{\overline \ell}
\newcommand\kb{\bar k}

\newcommand\ba{\bar \alpha}
\newcommand\bb{\bar \beta}
\newcommand\bc{\bar \gamma}
\newcommand\bq{\bar q}
\newcommand\bW{\bar W}
\newcommand\hW{\hat W}
\newcommand\bv{\bar y}
\newcommand\bB{\bar B}
\newcommand\bU{\bar U}
\newcommand\bT{\bar T}

\newcommand\hl{\hat \ell}
\newcommand\hL{\hat E}
\newcommand\hE{\hat y}
\newcommand\hF{\hat z}
\newcommand\hz{\hat \zeta}
\newcommand\cB{\mathcal{B}}
\newcommand\cC{\mathcal{C}}
\newcommand\cD{\mathcal{D}}
\newcommand\cG{\mathcal{G}}

\newcommand\eq[1]{\stepcounter{thm} \renewcommand{\theequation}{\thethm}
 \begin{equation}#1\end{equation}}

\newcommand\textmatrix[1]{\bigl(\begin{smallmatrix}#1\end{smallmatrix}\bigr)}

\begin{document}

\title[Exceptional Polynomials] {A New Family of Exceptional 
Polynomials in Characteristic Two}

\author{Robert M. Guralnick}
\address{Department of Mathematics, University of Southern California,
Los Angeles, CA 90089--2532, USA}
\email{guralnic@usc.edu}

\author{Joel E. Rosenberg}
\address{
  Center for Communications Research,
  4320 Westerra Court,
  San Diego, CA 92121--1967, USA
}
\email{joelr@ccrwest.org}

\author{Michael E. Zieve}
\address{
Center for Communications Research,
805 Bunn Drive,
Princeton, NJ 08540--1966, USA
}
\email{zieve@math.rutgers.edu}
\urladdr{http://www.math.rutgers.edu/$\sim$zieve/}

\thanks{We thank the referee for useful advice on notation.
The first author was partially supported by NSF grant
DMS 0653873}

\begin{abstract}
We produce a new family of polynomials $f(X)$ over fields $k$ of
characteristic $2$ which are exceptional, in the sense that
$f(X)-f(Y)$ has no absolutely irreducible factors in $k[X,Y]$
except for scalar multiples of $X-Y$; when $k$ is finite, this
condition is equivalent to saying that the map $\alpha\mapsto f(\alpha)$
induces a bijection on an infinite algebraic extension of $k$.
Our polynomials have degree $2^{e-1}(2^e-1)$,
where $e>1$ is odd.  We also prove that this completes the classification
of indecomposable exceptional polynomials of degree not a power of
the characteristic.
\end{abstract}

\maketitle


\section{Introduction}

Let $k$ be a field of characteristic $p\ge 0$, let
$f(X)\in k[X]\setminus k$, and let $\kb$ be an algebraic closure of $k$.
A polynomial in $k[X,Y]$ is called {\em absolutely irreducible}
if it is irreducible in $\kb[X,Y]$.
We say $f$ is \textit{exceptional} if $f(X)-f(Y)$ has no absolutely
irreducible factors in $k[X,Y]$ except for scalar multiples
of $X-Y$.  If $k$ is finite, this condition is
equivalent to saying that the map $\alpha\mapsto f(\alpha)$ induces a bijection on 
an infinite algebraic extension of $k$~\cite{Co,DL}.  Via this property, exceptional
polynomials have been used to construct remarkable examples of various types of objects:
curves whose Jacobians have real multiplication~\cite{TTV}, Galois extensions of number
fields with group $\PSL(q)$~\cite{DaM}, maximal curves over finite fields
\cite{CH,Oz}, families of character sums with small average value \cite{DL},
difference sets \cite{Dil,DD}, binary sequences with ideal
autocorrelation \cite{Dil}, almost perfect nonlinear power functions
\cite{Dob,Dob2,Dil2}, bent functions \cite{X,DD}, and double-error correcting
codes \cite{Dil2}.

Trivially any linear polynomial is exceptional.
The simplest nontrivial examples are the multiplicative
polynomials $X^d$ (which are exceptional when $k$ contains no $d$-th roots
of unity except $1$) and the additive polynomials $\sum \alpha_i X^{p^i}$
(which are exceptional when they have no nonzero root in $k$).
Dickson~\cite{Di} showed that certain variants of these polynomials are
also exceptional in some situations: the Dickson polynomials $D_d(X,\alpha)$
(with $\alpha\in k$), which are defined by $D_d(Y+\alpha/Y,\alpha)=Y^d+(\alpha/Y)^d$;
and the subadditive polynomials $S(X)$, which satisfy $S(X^m)=L(X)^m$
with $L$ an additive polynomial and $m$ a positive integer.  For nearly 100 years,
the only known exceptional polynomials were compositions of these
classical examples.

Klyachko~\cite{Kl}
showed that compositions of these polynomials yield all exceptional polynomials
of degree not divisible by $p$, and also all exceptional polynomials of
degree $p$.  A vast generalization of this result was proved by
Fried, Guralnick and Saxl~\cite{FGS}, which greatly restricted the possibilities
for the \emph{monodromy groups} of
exceptional polynomials.  We recall the relevant terminology: let $x$ be
transcendental over $k$.  We say 
$f(X)\in k[X]\setminus k$ is \emph{separable} if the field extension
$k(x)/k(f(x))$
is separable, or equivalently $f'(X)\ne 0$.  For a separable $f(X)\in k[X]$, 
let $E$ be the Galois closure of $k(x)/k(f(x))$.  The arithmetic 
monodromy group of $f$ (over $k$) is $\Gal(E/k(f(x)))$; the geometric
monodromy group of $f$ is $\Gal(E/\ell(f(x)))$, where $\ell$ is the 
algebraic closure of $k$ in $E$.  If $k$ is finite, then the composition
$b\circ c$ of two polynomials $b,c\in k[X]$ is exceptional if and only if 
both $b$ and $c$ are exceptional~\cite{DL}.  Thus, the study of 
exceptional polynomials over finite fields reduces to the case of
\emph{indecomposable} polynomials, i.e., polynomials which are not 
compositions of lower-degree polynomials.
For extensions of these results to infinite fields and to maps between
other varieties, see \cite{GS,GTZ,LMZ,Ma}.
Fried, Guralnick and Saxl proved the following result about the monodromy
groups of an indecomposable exceptional polynomial~\cite{FGS,GS}:

\renewcommand{\theenumi}{\roman{enumi}}
\renewcommand{\labelenumi}{(\theenumi)}
\begin{thm}\label{fgsthm}
Let $k$ be a field of characteristic $p$, and let $f(X)\in k[X]$
be separable, indecomposable, and exceptional of degree $d>1$.  Let $A$
be the arithmetic monodromy group of $f$.  Then one of the following 
holds.
\begin{enumerate}
\item $d\ne p$ is prime, and $A$ is solvable.
\item $d=p^e$ and $A$ has a normal elementary abelian subgroup $V$
of order $p^e$.
\item $p\in\{2,3\}$, $d=p^e(p^e-1)/2$ with $e>1$ odd, and
$A\cong\PGammaL(p^e)=\PGL(p^e)\rtimes\Gal(\F_{p^e}/\F_p)$.
\end{enumerate}
\end{thm}

It remains to determine the polynomials corresponding to these
group theoretic possibilities.  Case (i) is completely understood:
up to compositions with linear
polynomials, one just gets the Dickson polynomials $D_d(X,\alpha)$
(see~\cite[Appendix]{Mu-Schur} or \cite{Kl}).
In case (ii), we have $G=VG_1$ for some $G_1$; this case includes
the additive polynomials (where $G_1=1$) and
the subadditive polynomials (where $G_1$ is cyclic).  In joint work with
M\"uller \cite{GM,GMZ}, we have found families of case (ii) examples
in which $G_1$ is dihedral \cite{GM,GMZ}.
Moreover, in all known examples in case (ii), the fixed field $E^V$
has genus zero; conversely, we show in \cite{GMZ} that there
are no further examples in which $E^V$ has genus zero or one.
We suspect there are no other examples in case (ii):
for if $E^V$
has genus $g>1$ then $G_1$ will be a group of automorphisms of $E^V$
whose order is large compared to $g$, and there are not many possibilities
for such a field $E^V$.  We hope to complete the analysis of case (ii)
in a subsequent paper.
The present paper addresses case (iii).

In the two years following \cite{FGS}, examples were found in case (iii)
for each $p\in\{2,3\}$ and each odd $e>1$ \cite{CM,LZ,Mu}.
In the companion paper~\cite{GZ}, we show that twists of these examples
comprise all examples in case (iii), except possibly in the following
situation: $p=2$, $G=\SL(2^e)$, and the extension $k(x)/k(f(x))$
is wildly ramified over at least two places of $k(f(x))$.
In the present paper we conclude the treatment of case
(iii) by handling this final ramification setup.  In 
particular, we find a new family of exceptional polynomials.  Our main 
result is the following, in which we say polynomials $b,c\in k[X]$ are 
$k$-\emph{equivalent} if there are linear polynomials $\ell_1,\ell_2\in k[X]$ such
that $b=\ell_1\circ c\circ \ell_2$:

\begin{thm}
\label{thm-intro}
Let $k$ be a field of characteristic $2$. 
Let $q=2^e > 2$.  For $\alpha\in k\setminus\F_2$, define
\[
f_{\alpha}(X):=\left(\frac{\T(X)+\alpha}{X}\right)^q\cdot\left(\T(X)+\frac{\T(X)+\alpha}{\alpha+1}
\cdot \T\Bigl(\frac{X(\alpha^2+\alpha)}{(\T(X)+\alpha)^2}\Bigr)\right),
\]
where $\T(X)=X^{q/2}+X^{q/4}+\dots+X$.  Then the map $\alpha\mapsto f_{\alpha}$
defines a bijection from $k\setminus\F_2$ to the set of $k$-equivalence classes
of separable polynomials $f \in k[X]$ of degree $q(q-1)/2$ satisfying
\begin{enumerate}
\item the geometric monodromy group of $f$ is\/ $\SL(q)$; and
\item the extension $k(x)/k(f(x))$ is wildly ramified over at least two
places of $k(f(x))$.
\end{enumerate}
Every $f_{\alpha}$ is indecomposable.  Moreover, $f_{\alpha}$ is
exceptional if and only if $e$ is odd and
$k \cap \Fq = \F_2$.
\end{thm}

The strategy of our proof is to identify the curve $\cC$ corresponding to
the Galois closure $E$ of $k(x)/k(f(x))$, for $f$ a polynomial satisfying
(i) and (ii).  It turns out that $\cC$ is geometrically
isomorphic to the smooth plane curve $y^{q+1}+z^{q+1}=\T(yz)+\alpha$.

A key step in our proof is the computation
of the automorphism groups of curves of the form $v^q + v = h(w)$,
with $h$ varying over a two-parameter family of rational functions.  Our
method for this computation is rather general, and applies to many families
of rational functions~$h$.

As noted above, Theorem~\ref{thm-intro} completes the classification of
non-affine indecomposable exceptional polynomials:

\begin{cor}
\label{cor-intro}
Let $k$ be a field of characteristic $p\ge 0$.  Up to
$k$-equivalence, the separable indecomposable exceptional polynomials over $k$
which lie in cases (i) or (iii) of Theorem~\ref{fgsthm} are precisely:
\begin{enumerate}
\item for any $p$, the polynomial $X^d$ where $d\ne p$ is prime and
$k$ contains no $d$-th roots of unity except $1$;
\item for any $p$, the polynomial
\[
D_d(X,\alpha) := \sum_{i=0}^{\lfloor d/2\rfloor} \frac{d}{d-i}\binom{d-i}{i}
(-\alpha)^i X^{d-2i}
\]
where $d\ne p$ is prime,
$\alpha\in k^*$, and $k$ contains no elements of the form
$\zeta+1/\zeta$ with $\zeta$ being a primitive $d$-th root of unity in
$\bar{k}$;
\item for $p=2$ and $q=2^e>2$ with $e$ odd and $k\cap\F_q=\F_2$,
the polynomial $f_{\alpha}(X)$ where $\alpha\in k\setminus\F_2$;
\item for $p=2$ and $q=2^e>2$ with $e$ odd and $k\cap\F_q=\F_2$,
 the polynomial
\[
X\Bigl(\sum_{i=0}^{e-1}(\alpha X^n)^{2^i-1}\Bigr)^{(q+1)/n}
\]
where $n$ divides $q+1$ and $\alpha\in k^*$;
\item for $p=3$ and $q=3^e>3$ with $e$ odd and $k\cap\F_q=\F_3$,
the polynomial
\[
X(X^{2n}-\alpha)^{(q+1)/(4n)}\left(\frac{(X^{2n}-\alpha)^{(q-1)/2}
+\alpha^{(q-1)/2}}{X^{2n}}\right)^{(q+1)/(2n)}
\]
where $n$ divides $(q+1)/4$ and $\alpha\in k^*$ has image in
$k^*/(k^*)^{2n}$ of even order.
\end{enumerate}
\end{cor}

The contents of this paper are as follows.
In the next section we prove some useful
results about ramification groups.  In Section~\ref{sec-GZ} we record results
from \cite{GZ} which describe the ramification (including the higher 
ramification groups) in $E/k(f(x))$.  In Section~\ref{sec-B-curves}
we classify curves which admit $B$, the group of upper triangular matrices
in $\SL(q)$, as a group of automorphisms with our desired ramification
configuration.  There is a two-parameter family of such curves.  In 
Sections~\ref{sec-aut-groups} and \ref{sec-G-curves} we
determine the automorphism groups of the curves in this
family (which turn out to be either $B$ or $\SL(q)$).  The curves
with automorphism group $\SL(q)$ form a one-parameter subfamily.
The group theoretic data yields the existence and uniqueness
of the desired polynomials.  In particular, it shows we cannot have
$k=\F_2$; in Section~\ref{sec-nonexistence} we give a different, more
direct proof of this fact.  In the final two sections we consider different
forms of the curves, and in particular we determine a smooth plane model.
We then use this model to explicitly compute the polynomials, and we conclude
the paper by proving Theorem~\ref{thm-intro} and Corollary~\ref{cor-intro}.

\textbf{Notation.}
Throughout this paper, all curves are assumed to be smooth, projective, and
geometrically irreducible.  We often define a curve by giving an affine
plane model, in which case we mean the completion of the normalization
of the stated model.
Also, in this case we describe points on the curve by giving the 
corresponding points on the plane model.

A \emph{cover} is a separable nonconstant morphism between curves.
If $\rho\colon \cC\to \cD$ is a cover of curves over a field $k$, by a `branch 
point' of $\rho$ we mean a point of $\cD(\kb)$ which is ramified in 
$\rho\times_k \kb$ (for $\kb$ an algebraic closure of $k$).  In 
particular, branch points need not be defined over $k$, but
the set of branch points is preserved by the absolute Galois group of $k$.  If $f\in k[X]$
is a separable polynomial, then we refer to the branch points of the corresponding cover
$f\colon\Pone\to\Pone$ as the branch points of $f$.

If $\rho\colon \cC\to \cD$ is a Galois cover, and $P$ is a point of $\cC$, then the
ramification groups at $P$ (in the lower
numbering, as in~\cite{Se}) are denoted $I_0(P), I_1(P), \ldots$,
or simply $I_0, I_1, \ldots$.
Here $I_0$ is the inertia group and $I_1$ is its Sylow $p$-subgroup.
We refer to $I_1$ as the first ramification group, $I_2$ as the second, and so on.

We reserve the letter $x$ for an element transcendental over the field $k$.
Throughout this paper we write $q=2^e$ where $e>1$.
We use the following notation for subgroups of $\SL(q)$.
The group of diagonal matrices is denoted $T$.
The group of upper triangular matrices is denoted $B$.  
The group of elements of $B$ with 1's on the diagonal
is denoted $U$.  The two-element group generated by 
$\textmatrix{1&1\\0&1}$ is denoted $W$.  Finally, $\T(X)$ denotes the
polynomial $X^{q/2}+X^{q/4}+\dots+X$.


\section{Ramification in Galois $p$-power covers}
\label{sec-ordinary}

In this section we prove a useful result (Corollary~\ref{cor-ordinary})
about ramification groups in Galois covers of degree a power of the
characteristic.  We give two proofs, each of which provides additional
information.  Throughout this section, $\bl$ is an algebraically closed
field of characteristic $p>0$.

\begin{prop}
\label{prop-ramification}
Let $\rho\colon \cC\to \cD$ and $\rho'\colon \cC\to \cB$ be Galois covers
of curves over $\bl$.  Let $n$ and $r$ be positive integers.
Suppose that $\cB\cong\Pone$ and the degree of $\rho'$
is $p^r$.  If all $n$-th ramification groups of $\rho'$ are
trivial, then the same is true of $\rho$.
\end{prop}

\begin{proof}
Since $\Pone$ has no nontrivial unramified covers, and any
ramified Galois cover of $p$-power degree has a nontrivial
first ramification group, the hypotheses imply $n\geq 2$.  
Without loss, we may assume $\rho$ has degree $p$.
First assume $\cC$ has genus greater than $1$, so that
$\Aut \cC$ is finite.  Let $H$ be a Sylow $p$-subgroup of
$\Aut \cC$ which contains $\Gal(\rho')$.  By replacing $\rho$
by one of its $(\Aut \cC)$-conjugates, we may assume that
$H$ contains $\Gal(\rho)$ as well.  Then the cover 
$\cB\to \cC/H$ induced from $\cC\to \cC/H$ and $\cC\to \cB$ is the 
composition of a sequence of Galois degree-$p$
covers $\cB=\cB_0\to \cB_1\to\dots\to \cB_m=\cC/H$.  Since $\cB\cong\Pone$,
each $\cB_i\to \cB_{i+1}$ has trivial second ramification groups
(by Riemann-Hurwitz).  Thus, each $n$-th ramification group of 
$\cC\to \cB_{i+1}$ is also an $n$-th ramification group of
$\cC\to \cB_i$ (by \cite[Prop.~IV.3]{Se}).  By induction,
$\cC\to \cB_m=\cC/H$ has trivial $n$-th
ramification groups, whence the same is true of $\cC\to \cD$.

If $\cC$ has genus $0$, then $\rho$ is a degree-$p$ cover between
genus-$0$ curves and hence has trivial second ramification groups.

Finally, assume $\cC$ has genus $1$.  Pick a point of $\cC$ with
nontrivial inertia group under $\rho'$, and let $J$ be an order-$p$
subgroup of this inertia group.  Then $\cC/J\cong\Pone$, so by replacing
$\cB$ by $\cC/J$ we may assume $\rho'$ has degree $p$.  If $p>3$ then
no such $\rho'$ exists (e.g., by Riemann-Hurwitz).  If $p=3$ then
any Galois degree-$p$ map $\cC\to \cD$ is either unramified (with $\cD$ of
genus $1$) or has a unique branch point (with $I_2\ne I_3$ and $\cD$ of genus $0$).
Henceforth assume $p=2$.  Then a degree-$p$
map $\cC\to \cD$ is either unramified (with $\cD$ of genus $1$) or has
precisely two branch points (each with $I_1\ne I_2$ and $\cD$ of genus $0$) or has
a unique branch point (with $I_3\ne I_4$ and $\cD$ of genus $0$).
If there is a unique branch point then $\cC$ is isomorphic to the curve
$y^2+y=z^3$.  Since the corresponding elliptic curve has trivial $2$-torsion,
it follows that a degree-2 function on this curve cannot have two branch points.
This completes the proof.
\end{proof}

We will use the following result, which follows from
Proposition~\ref{prop-ramification} and standard results about
ramification groups (cf.~\cite[\S IV2]{Se}).

\begin{cor}
\label{cor-ordinary}
Let $\rho\colon \cC\to \cD$ and $\rho'\colon \cC\to\Pone$ be Galois covers
of curves over $\bl$.  Suppose that $\rho'$ has degree a power of $p$,
and that all second ramification groups of $\rho'$ are trivial.
If $I_1$ and $I_2$ are the first and second ramification groups of $\rho$
at some point of $\cC$, then
$I_1$ is elementary abelian, $I_2=1$, and $I_1$ is its own
centralizer in $I$.
\end{cor}

We now give a different proof of this corollary, which generalizes
the corollary in a different direction than does 
Proposition~\ref{prop-ramification}.  For a curve $\cC$ over $\bl$,
let $p_\cC$ denote the $p$-rank of $\cC$ (i.e., the rank of the $p$-torsion 
subgroup of the Jacobian of $\cC$).  Let $g_\cC$ denote the 
genus of $\cC$.  These quantities are related by $p_\cC \le g_\cC$.
Recall that $\cC$ is called {\it ordinary} if $p_\cC = g_\cC$.  
We first record a standard basic fact.

\begin{lemma} 
Let $\theta\colon \cC\to \cD$ be a cover of curves over $\bl$.
If $\cC$ is ordinary, then $\cD$ is ordinary.
\end{lemma}

This lemma and the next one are proved in \cite[Thm.\ 1.2]{Pi}.
The strategy for proving the next lemma comes from
\cite[Thm.\ 2]{Na}.

\begin{lemma}
\label{lemma-ordinary}
Let $\theta\colon \cC \to \cD$ be a Galois cover (of curves over $\bl$) 
whose Galois group $H$ is a $p$-group.  Then $\cC$ is ordinary if and
only if both
\begin{enumerate}
\item $\cD$ is ordinary; and
\item every branch point of $\theta$ has trivial second ramification group.
\end{enumerate}
\end{lemma}

\begin{proof}
We use the Deuring-Shafarevich formula (\cite[Thm.~4.2]{Su}):
\begin{equation*} 
\frac{p_\cC - 1}{|V|} = p_\cD - 1 + \sum_{Q \in \cD} \Bigl(1 -
 \frac{1}{e_Q}\Bigr), 
\end{equation*}
where $e_Q$ is the ramification index of $\theta$ at the point $Q$.

The Riemann-Hurwitz formula yields
\begin{equation*} 
\frac{g_\cC - 1}{|V|} = g_\cD - 1 + \sum_{Q \in \cD} \Bigl(1 -
\frac{1}{e_Q}\Bigr) + s,
\end{equation*}
where $s$ is the contribution from the second and higher ramification
groups.  Note that $s\ge 0$, with equality if and only if all
second ramification groups are trivial.

Since $p_\cD\le g_\cD$, we conclude that $p_\cC=g_\cC$ holds if and only if 
$p_\cD=g_\cD$ and $s=0$.
\end{proof}

\begin{proof}[Alternate proof of Corollary~\ref{cor-ordinary}]
By applying the previous result with $\theta=\rho'$, we see that
$\cC$ is ordinary.  Applying it with $\theta=\rho$ shows that $I_2=1$,
and then the remaining assertions follow from standard properties
of the higher ramification groups.
\end{proof}

\begin{remark}
It would be interesting to refine the above alternate proof of
Corollary \ref{cor-ordinary} to prove Proposition \ref{prop-ramification}.
Such a refinement would likely require a refinement of the
Deuring-Shafarevich formula that involves finer invariants than just
the $p$-rank.  However, we do not know such a refined formula.
We thank Hendrik Lenstra for suggesting this possibility.
\end{remark}


\section{Previous results}
\label{sec-GZ}

We will use the following result from the companion paper \cite[Lemma 2.7]{GZ}.
Recall our convention that $x$ is
transcendental over $k$; also $B$ is the group of upper triangular matrices in
$\SL(q)$, and $W=B\cap\SL(2)$.

\begin{lemma}
\label{lemma-GZ}
Let $k$ be a perfect field of characteristic $2$, and let $q=2^e$ with $e>1$.
Suppose $f\in k[X]$ is a separable polynomial
of degree $q(q-1)/2$ which satisfies conditions \emph{(i)} and \emph{(ii)} of
Theorem~\ref{thm-intro}.  Let $E$ be the Galois closure of
$k(x)/k(f(x))$, and let $\ell$ be the algebraic closure of
$k$ in $E$.
Then $E/\ell(f(x))$ has precisely two ramified places, both of degree one,
and the corresponding inertia groups are $B$ and $W$ (up to conjugacy).
Moreover, the second ramification group over each ramified place
is trivial, and $f$ is indecomposable.  The degree $[\ell\col k]$ divides $e$,
and $f$ is exceptional
if and only if $e$ is odd and $[\ell\col k]=e$.  Finally, there is a 
curve $\cC_0$ over $k$ such that $\ell.k(\cC_0)\cong_{\ell} E$.
\end{lemma}

The following consequence of Lemma~\ref{lemma-GZ} describes the
ramification in $\cC\to \cC/B$, where $\cC=\cC_0\times_k \ell$.
This too was proved in the companion paper \cite[Cor.~2.8]{GZ}.
Here $T$ is the group of diagonal matrices in $\SL(q)$.

\begin{cor}
\label{cor-B}
If $\cC$ is a curve over $\ell$ for which $\ell(\cC)=E$, then the following hold:
\begin{enumerate}
\item $B$ acts as a group of $\ell$-automorphisms on $\cC$;
\item the quotient curve $\cC/B$ has genus zero;
\item the cover $\cC \to \cC/B$ has exactly three branch points;
\item the inertia groups over these branch points are $B$, $T$, and
      $W$ (up to conjugacy); and
\item all second ramification groups in the cover $\cC \to \cC/B$ are 
      trivial.
\end{enumerate}
\end{cor}

We now record some standard facts about subgroups of $\SL(q)$;
see for instance \cite[\S260]{Di2}, \cite[\S3.6]{Sz}, or \cite[App.]{GZ}.
Here $U$ is the group of elements of $B$ whose diagonal entries are $1$.

\begin{lemma}
\label{lemma-gt}
$B=U\rtimes T$ is the semidirect product of the normal subgroup $U$
by the cyclic subgroup $T$.  All involutions in $B$ are conjugate.
All subgroups of $B$ of order $q-1$ are conjugate.
For $j\in\{1,-1\}$,
all subgroups of\/ $\SL(q)$ of order $2(q+j)$ are conjugate, and these
subgroups are dihedral and are maximal proper subgroups of\/ $\SL(q)$.
The normalizer of $W$ in\/ $\SL(q)$ is $U$.
There is no group strictly between $B$ and\/ $\SL(q)$.
\end{lemma}


\section{$B$-curves}
\label{sec-B-curves}

Let $\ell$ be a perfect field of characteristic $2$.
In this section we describe the curves $\cC$ over $\ell$ which
admit a $B$-action as in Corollary~\ref{cor-B}.  We will show that the only
such curves $\cC$ are the curves $\cC_{\alpha,\beta}$ defined as follows.
For any $\alpha,\beta\in\ell^*$, let $\cC_{\alpha,\beta}$ be the curve defined by

\eq{
\label{eq-Cab}
v^q + v = (\alpha + \beta)w + w^q \T\Bigl(\frac{\beta}{1+w^{q-1}}\Bigr),}
where\/ $\T(X):=X^{q/2}+X^{q/4}+\dots+X$.
Note that $\cC_{\alpha,\beta}$ is geometrically irreducible, since the 
left side of (\ref{eq-Cab}) is a polynomial in $v$ and the right 
side is a rational function in $w$ with a simple pole 
(at $w=\infty$, with residue $\alpha$).

\begin{thm}
\label{thm-B-curves}
Suppose $\cC$ is a curve over $\ell$ satisfying the five properties in 
Corollary~\ref{cor-B}.  Then $\ell\supseteq\Fq$ and $\cC\cong \cC_{\alpha,\beta}$
for some $\alpha,\beta\in \ell^*$.
\end{thm}

Indeed, suppose $\cC$ satisfies the properties of Corollary~\ref{cor-B}.
Since the inertia groups $B$, $T$, and $W$ are not conjugate, the 
corresponding branch points are $\ell$-rational,
so for a suitably chosen coordinate $t$ on $\cC/B$ they are $\infty$, $0$, and $1$,
respectively.  Note that $\cC/U \to \cC/B$ is a cyclic cover of degree $q-1$
which is totally ramified over $\infty$ and $0$, and unramified elsewhere.
By Riemann-Hurwitz, $\cC/U$ has genus zero.
Each of the $q-1$ order-$2$ subgroups of $U$ is conjugate to $W$,
and is thus an inertia group in $\cC \to \cC/B$, hence also in $\cC \to \cC/U$.
Thus there must be at least $q-1$ distinct places of $\cC/U$ lying over
the place $t=1$ of $\cC/B$, so all of these places must be rational.
Choose a coordinate $w$ on $\cC/U$ such that, in the cover $\cC/U \to \cC/B$, 
the points $\infty$, $0$, and $1$ map to $\infty$, $0$, and $1$, respectively.
Then $\ell(\cC/U)=\ell(w)$ and $\ell(\cC/B)=\ell(t)$ where $t=w^{q-1}$.
Since $\cC/U \to \cC/B$ is Galois, $\ell$ contains $\Fq$.

In these coordinates, the branch points of the cover
$\cC \to \cC/U$ are $\infty$ and the $q-1$ elements of $\Fq^*$ (i.e., the 
points over $t=1$)
with the corresponding inertia groups being $U$ and its
$q-1$ subgroups of order~$2$.

Let $\cC_1=\cC/H$, where $H$ is a maximal subgroup of $U$.
Since $\cC\to \cC/U$ has no nontrivial second ramification groups, the same
is true of $\cC_1\to \cC/U$, so (since $\ell$ is perfect) $\cC_1$ is defined
by an equation of the form
\eq{
\label{veq}
y^2 + y = \alpha w + \sum_{\zeta \in \Fq^*}  
\frac{\beta_\zeta\zeta}{w + \zeta} + \gamma}
for some $y \in \ell(\cC)$ and $\alpha,\beta_\zeta,\gamma\in \ell$.
Note that $\alpha \ne 0$ (since $w=\infty$ is a branch point).
Clearly $\beta_\zeta\ne 0$ if and only if $w=\zeta$ is a branch point
of the cover $\cC_1 \to \cC/U$, and the latter holds if and only if $H$
does not contain the inertia group of $w=\zeta$ in $\cC \to \cC/U$.  Thus,
$\beta_{\zeta}$ is nonzero for precisely $q/2$ values $\zeta$.

Let $\Gamma$ be the set of elements $z\in \ell(\cC)$ for which
\begin{equation*}
z^2 + z = \ba(z)w + \sum_{\zeta \in \Fq^*} 
  \frac{\bb_\zeta(z) \zeta}{w + \zeta} + \bc(z)
\end{equation*}
with $\ba(z), \bb_\zeta(z), \bc(z)\in \ell$.
Note that $\ba(z)$, $\bb_\lambda(z)$, and $\bc(z)$ are uniquely
determined by
$z$, and each of them defines a homomorphism $\Gamma\to \ell$.
Let $\Gamma_0 = \Gamma \cap \ell(w)$; considering orders of poles,
we see that $\Gamma_0=\ell$.

Since $B=UT$, restriction to $\cC/U$ induces an isomorphism $T\cong B/U$,
so $T=\{\phi_\eta : \eta\in\F_q^*\}$ where $\phi_\eta(w)=\eta w$.
Clearly $\Gamma$ is $T$-invariant.
The following lemma enables us to choose $y$ so that $Ty\cup\{0\}$ is a group.

\begin{lemma}
There exists an order-$q$ subgroup\/ $\Gamma_1$ of\/ $\Gamma$ such that\/
$\Gamma=\Gamma_0 \oplus \Gamma_1$ and the nonzero elements of\/ $\Gamma_1$
comprise a single $T$-orbit.
\end{lemma}

\begin{proof}
The map $\theta\colon z+\Gamma_0\mapsto \ell(w,z)$
defines a surjective $T$-set homomorphism between $\Gamma/\Gamma_0\setminus\{0\}$
and the set $\Lambda$ of degree-2 extensions of $\ell(\cC/U)$ contained in $\ell(\cC)$.
We first prove injectivity of $\theta$:
suppose $z_1, z_2\in\Gamma\setminus\Gamma_0$ satisfy $\ell(w,z_1)=\ell(w,z_2)$.
Then the nonidentity element of $\Gal(\ell(w,z_1)/\ell(w))$ maps
$z_1\mapsto z_1+1$ and $z_2\mapsto z_2+1$, hence fixes
$z_1+z_2$, so $z_1+z_2\in\Gamma_0=\ell$.  Hence $\theta$ is injective.  Since $\Lambda$
is a transitive $T$-set of size $q-1$, it follows that $|\Gamma/\Gamma_0|=q$
and $T$ acts transitively on $\Gamma/\Gamma_0\setminus\{0\}$.  Finally, since
$|T|$ is odd and both $\Gamma$ and $\Gamma_0$ are $T$-invariant elementary
abelian $2$-groups, Maschke's theorem (\cite[12.9]{As}) implies there is a
$T$-invariant group $\Gamma_1$ such that $\Gamma=\Gamma_0\oplus\Gamma_1$,
and $|\Gamma_1|=|\Gamma/\Gamma_0|=q$.
\end{proof}

By replacing $y$ by $y+\delta$ for some $\delta\in\Gamma_0$, we
may assume that $y$ is in $\Gamma_1$.
Applying $\phi_\eta$ to (\ref{veq}), we see that
$y_\eta:=\phi_\eta(y)$ satisfies
\begin{equation*}
y_{\eta}^2 + y_{\eta} = \alpha\eta  w + 
\sum_{\zeta \in \Fq^*}  
\frac{\beta_{\zeta}\eta^{-1}\zeta}{w + \eta^{-1}\zeta}+\gamma.
\end{equation*}
Thus, $\bc(y_{\eta})=\gamma$ and $\ba(y_{\eta}) = \alpha \eta$ and
$\bb_{\eta^{-1}\zeta}(y_{\eta})=\beta_{\zeta}$.
Since the homomorphism $z\mapsto \bc(z)$ is constant on the nonzero
elements of the group $\Gamma_1$, it follows that $\gamma=0$.

Since $\Gamma_1 = \{y_\eta\}\cup\{0\}$ is closed under
addition, $y_{\eta} + y_{\eta'} = y_{\eta''}$ for
some $\eta''$.  Comparing images under $\ba$ yields that 
\begin{equation*}
y_{\eta} + y_{\eta'} = y_{\eta + \eta'}.
\end{equation*}
Thus, 
\begin{equation*}
\beta_\zeta + \beta_\eta =
   \bb_1(y_{\zeta}) + \bb_1(y_{\eta}) = 
   \bb_1(y_{\zeta + \eta}) =
   \beta_{\zeta + \eta}.
\end{equation*}
Since $\beta_{\zeta}=0$ for exactly $q/2-1$ choices of $\zeta \in \Fq^*$,
this implies that $\beta_{\zeta}=0$ for $\zeta$ in some hyperplane
(i.e., index-$2$ subgroup) $\mathcal H$ of $\Fq$,
and $\beta_\zeta=\beta_{\zeta'}$ for $\zeta,\zeta'\not\in\mathcal H$.
Hence, $\bb_\zeta(y_\eta)=0$ for $\zeta\in\eta^{-1}\mathcal H$.
The hyperplanes $\eta^{-1}\mathcal H$ comprise all $q-1$ hyperplanes
in $\Fq$, so there is some $\eta$ for which $\eta^{-1}\mathcal H$ is
the set of roots of $\T(X):=X^{q/2}+X^{q/4}+\dots+X$.
Replacing $y$ by $y_\eta$, the equation for $\cC_1$ becomes
\begin{equation*}
y^2 + y = \alpha w + \beta \sum_{\zeta \in \Fq^*}  
\frac{\T(\zeta)\zeta}{w + \zeta}.
\end{equation*}
Note that $\alpha$ and $\beta$ are nonzero elements of $\ell$.

Since $\ell(\cC)$ is the Galois closure of $\ell(\cC_1)/\ell(w^{q-1})$,
it is uniquely determined by the choice of
$\alpha$ and $\beta$.  Thus, to conclude the proof of Theorem~\ref{thm-B-curves},
it suffices to show that (for each choice of $\alpha,\beta\in \ell^*$) the curve
$\cC_{\alpha,\beta}$ satisfies the hypotheses of the theorem, and that the
quotients of $\cC_{\alpha,\beta}$ by $B$ and by some order-$q/2$ subgroup induce
the above cover $\cC_1\to\Pone_{w^{q-1}}$.  The following lemma is clear:

\begin{lemma}
\label{lemma-BCab}
If $\ell$ contains\/ $\Fq(\alpha,\beta)$, then for any 
$\textmatrix{\gamma^{-1}&\delta\\0&\gamma}\in B$ there is a unique $\ell$-automorphism
of $\cC_{\alpha,\beta}$ mapping $w\mapsto \gamma^2 w$ and $v\mapsto \gamma^2 v + \gamma\delta$.
This correspondence defines an embedding
$B\hookrightarrow\Aut_\ell(\cC_{\alpha,\beta})$.
\end{lemma}

We now show that $\cC_{\alpha,\beta}$ (together with this action of $B$) has the
desired properties.

\begin{lemma}
\label{lemma-Cab}
The curve $\cC:=\cC_{\alpha,\beta}$ has genus $q(q-1)/2$.  
Moreover, the fixed fields $\ell(\cC)^U$ and $\ell(\cC)^B$ equal $\ell(w)$ and $\ell(w^{q-1})$,
and the cover $\cC \to \cC/B$ has precisely three branch points.
The inertia groups over these points are (up to conjugacy)
$B$, $T$, and $W$.  Also, the second ramification groups at all three
points are trivial.  Finally, if
$H=\{\textmatrix{1&\delta\\0&1} : \T(\delta)=0\}$,
then $\ell(\cC)^H = \ell(w,y)$ where 
\eq{\label{eq-v}
y^2 + y = \alpha w + \beta \sum_{\zeta \in \Fq^*} 
\frac{\T(\zeta)\zeta}{w - \zeta}.}
\end{lemma}

\begin{proof} 
It is clear that $\ell(\cC)^U=\ell(w)$ and $\ell(\cC)^B=\ell(t)$, where
$t:=w^{q-1}$.  Also, both $w$ and $\bv:=\T(v)$ are fixed by $H$, and
a straightforward calculation yields
\begin{equation*}
\bv^2 + \bv = (\alpha+\beta)w + w^q\T\Bigl(\frac{\beta}{t+1}\Bigr) = 
\alpha w  + \beta \sum_{\zeta \in \Fq^*}  
\frac{\T(\zeta)\zeta}{w - \zeta} + h + h^2
\end{equation*}
for an appropriate $h \in \ell(w)$.  Thus, $H$ fixes $w$ and 
$y:=\bv+h$, and $y$ satisfies (\ref{eq-v}).  Since $y\notin\ell(w)$,
it follows that $\ell(\cC)^H=\ell(w,y)$.  Note that the genus of 
$\ell(w,y)$ is $q/2$, since the right hand side of (\ref{eq-v}) has 
precisely $1+q/2$ poles and they are all simple.

Let $\cD=\cC/U$ and $\cB=\cC/B$, so $\ell(\cD)=\ell(w)$ and $\ell(\cB)=\ell(t)$.
The cover $\cD\to \cB$ is only ramified at $w=0$ and $w=\infty$, and
is totally ramified at both of these points.  The cover $\cC\to \cD$
can only be ramified at points with $v=\infty$, hence at points with
$w\in\Fq^*$ or $w=\infty$.  The point $w=\infty$ of $\cD$ is totally
ramified in $\cC\to \cD$, since $w$ is a simple pole of the right hand
side of (\ref{eq-Cab}).  The points $w\in\Fq^*$ of $\cD$ all lie over
the point $t=1$ of $\cB$, and precisely $q/2$ of these points are ramified
in $\cC/H\to \cD$.  Since $T$ permutes transitively both the $q-1$ points in
$\cD$ over $t=1$ and the $q-1$ index-2 subgroups of $U$, we see that each
such point ramifies in precisely $q/2$ of the covers 
$\cC/V \to \cD$ as $V$ ranges over the $q-1$ index-2 subgroups of $U$.
This implies that each $w\in\Fq^*$ has ramification index $2$ in
$\cC\to \cD$.  Thus, the only branch points of the cover $\cC\to \cC/B$ are
$\infty$, $0$, and $1$, and the corresponding ramification indices are
$q(q-1)$, $q-1$, and $2$.  Hence, up to conjugacy, the corresponding
inertia groups are $B$, $T$, and $W$.  Moreover, since the second
ramification groups of $\cC/H\to \cD$ are trivial, the same is true of
every $\cC/V\to \cD$, and hence of $\cC\to \cD$.  It follows from
Riemann-Hurwitz that $\cC$ has genus $q(q-1)/2$.
\end{proof}

\noindent
This concludes the proof of Theorem~\ref{thm-B-curves}.


\section{Automorphism groups of $B$-curves}
\label{sec-aut-groups}

Let $\bl$ be an algebraically closed field of characteristic $2$,
let $\alpha,\beta\in\bl^*$, and put $\cC:=\cC_{\alpha,\beta}$ as in (\ref{eq-Cab}).
By Lemma~\ref{lemma-Cab},
$\cC$ admits an action of $B$ satisfying the five properties of
Corollary~\ref{cor-B}.  In this section we prove that the automorphism
group of $\cC$ is either $B$ or $\SL(q)$.

Let $P_1$, $P_2$, and $P_3$ be points of $\cC$ whose stabilizers (in $B$)
are $B$, $W$, and $T$, respectively.
Let $\cG$ be the automorphism group of $\cC$.

\begin{lemma}
\label{lemma-B-curves}
Let $V \le U$ be a subgroup with $|V|> 2$.
Then $N_{\cG}(V)\le B$ and $|\cG\col B|$ is odd.  Moreover, $B$ is the stabilizer
of $P_1$ in $\cG$.  
\end{lemma}

\begin{proof} 
Let $\bB$ be the stabilizer of $P_1$ in $\cG$, and let $\bU$ be the
Sylow $2$-subgroup of $\bB$.  Corollary~\ref{cor-ordinary} implies that
$\bU$ is elementary abelian and that $\cC\to \cC/\bU$ has trivial second
ramification groups.

Write $|\bU|=q\bq$.  Since the
$q-1$ order-$2$ subgroups of $U$ are all conjugate under $B$, they
are all inertia groups in $\cC\to \cC/U$.  These subgroups are nonconjugate in
the abelian group $\bU$, so $\cC\to \cC/\bU$ has at least $q-1$ distinct
branch points not lying under $P_1$.
By Riemann-Hurwitz,
$2(q\bq + q(q-1)/2 - 1) = \sum_Q \ind(Q)$ where $Q$ varies over
the branch points of $\cC\to \cC/\bU$ and $\ind(Q)$ is the sum of the
different exponents (in the cover $\cC\to \cC/\bU$) of the points over
$Q$.  If $Q$ lies under $P_1$, then $Q$ is totally ramified so
$\ind(Q) = 2(q\bq-1)$.  Any branch point satisfies $\ind(Q) \ge q\bq$
(since $\cC\to \cC/\bU$ is a Galois cover with Galois group a $2$-group).  
Thus,
\begin{equation*}
2\Bigl(q\bq + \frac{q(q-1)}{2} - 1\Bigr) \ge  2(q\bq-1) + (q-1)q\bq,
\end{equation*}
or $q(q-1) \ge (q-1)q\bq$.  Hence $\bq=1$, so $\bU=U$.

By Corollary~\ref{cor-ordinary}, $U$ is its own centralizer in $\bB$,
so conjugation induces a faithful action of $\bB/U$ on $U$ and thus
also on $U\setminus\{0\}$.
Since $\bB/U$ is cyclic, it follows that $|\bB/U|\le |U\setminus\{0\}|=
|B/U|$, so $\bB=B$.

Since we know the inertia groups of $\cC\to \cC/B$, we see that $P_1$
is the only point of $\cC$ fixed by $V$.  Thus, $N_{\cG}(V)$ fixes $P_1$,
so $N_{\cG}(V)\le B$.  In particular, $N_{\cG}(U)=N_B(U)=B$.
If $U$ is not a full Sylow $2$-subgroup of $\cG$, then (since 
$2$-groups are nilpotent) $|N_{\cG}(U)\col U|$ is even, a contradiction.
Thus, $|\cG\col B|$ is odd.
\end{proof}

\begin{lemma}
\label{MBequiv}
The following are equivalent:
\begin{enumerate}  
\item $\cG=B$;
\item $P_1$ and $P_3$ are in distinct $\cG$-orbits; 
\item $T=N_{\cG}(T)$; and
\item $|N_{\cG}(T)\col T| \ne 2$.
\end{enumerate}
\end{lemma}

\begin{proof}
By Lemma \ref{lemma-B-curves}, the intersection of the stabilizers (in $\cG$)
of $P_1$ and $P_3$ is $T$.  
Since distinct $B$-conjugates of $T$ intersect trivially, any nontrivial element of
$T$ fixes precisely two points of $\cC$ (namely $P_1$ and $P_3$).
Thus, $N_{\cG}(T)$ preserves $\Lambda:=\{P_1,P_3\}$, so either it acts transitively
on $\Lambda$ (and $|N_{\cG}(T)\col T|=2$) or else $N_{\cG}(T)=T$.  Hence conditions (iii)
and (iv) are equivalent, and they both follow from (ii).
If $\nu P_3 = P_1$ with $\nu\in \cG$, then
$T^{\nu}$ is contained in $B$, so $T^{\nu}=T^{\mu}$ for some $\mu\in B$;
but then $\mu^{-1}\nu\in N_{\cG}(T)\setminus T$.  Hence (ii) and (iii) are equivalent.

Clearly if $\cG=B$, all the remaining conditions are true.
So we assume the last three conditions and show that $\cG=B$.

Suppose $P_1$ and $P_3$ are in distinct $\cG$-orbits.  Let $I$ be
the stabilizer of $P_3$ in $\cG$, so $I=V\bT$ where $V$ is a normal
2-subgroup and $\bT$ is a cyclic group of odd order.  Since $I$ contains
$T$, by Schur-Zassenhaus $T$ is contained in an $I$-conjugate $\bT'$ of $\bT$,
so $I=V\bT'$.
Since $\bT'$ is cyclic, it normalizes $T$, so (by (iii)) $\bT'=T$.
Since $U$ is a Sylow $2$-subgroup of $\cG$ (Lemma \ref{lemma-B-curves}),
some conjugate $V'$ of $V$ is
contained in $U$; by our hypothesis on the inertia groups of
$\cC\to \cC/B$, either $|V|\le 2$ or $V'=U$.  But $V'\ne U$ because
$P_1$ and $P_3$ are in distinct $\cG$-orbits, and $|V|\ne 2$ since
$|I\col T|=2$ contradicts (iii).  Hence $I=T$.  Since any nontrivial element of $T$
fixes no point of $\cG P_3\setminus\{P_3\}$, it follows that $\cG$ acts on $\cG P_3$ as a
Frobenius group with Frobenius complement $T$; let $K$ be the Frobenius
kernel.  Since $K$ is a normal subgroup of $\cG$ that contains a Sylow
2-subgroup, $K$ contains every Sylow $2$-subgroup, so $U\le K$.  By
Lemma \ref{lemma-B-curves}, $N_{\cG}(U)=B$, so $N_K(U)=B\cap K=U$.
Nilpotence of the Frobenius kernel implies $K=U$, so $\cG=KT=B$.
\end{proof}

\begin{lemma}
\label{lemma-p2}
$W$ is the stabilizer of $P_2$ in $\cG$.
\end{lemma}

\begin{proof}
Let $\hW$ be the stabilizer of $P_2$ in $\cG$.
Let $\bW$ be the Sylow $2$-subgroup of $\hW$.
By Corollary \ref{cor-ordinary}, $\bW$ is elementary abelian and
is its own centralizer in $\hW$.  Thus, $\hW/\bW$ embeds in $\Aut(\bW)$.
If $\bW=W$, this implies that $\hW=W$.  Now assume that $\bW$ strictly
contains $W$; we will show that this leads to a contradiction.
Note that Lemma~\ref{MBequiv} implies $P_3\in \cG P_1$.

Let $C$ be the centralizer of $W$ in $\cG$.  Then $C$ contains $U$
and $\bW$, where $\bW\cap U=W$.  Let $\Lambda=CP_2$.
Since $W$ and $C$
commute, $W$ acts trivially on $\Lambda$, so $\Lambda\subseteq \{P_1\}\cup
UP_2$.  But $U$ is a Sylow $2$-subgroup of $\cG$, so it contains a conjugate
$\bW^{\nu}$ of $\bW$ in $\cG$, and since $|\bW|>2$ we must have $\nu P_2=P_1$.  Hence $\Lambda=
\{P_1\}\cup UP_2$.  The stabilizer of $P_1$ in $\cG$ is $B$, and the
stabilizer in $B$ of any element of $UP_2$ is $W$.  Thus any two-point
stabilizer of $C$ on $\Lambda$ is conjugate in $C$ to $W$, hence equals $W$,
so $C/W$ is a Frobenius group on $\Lambda$.
A Frobenius complement is $U/W$ (since $B\cap C=U$).
It is well known (and elementary in this case: cf.\ \cite[Thm.~3.4A]{DM})
that an abelian subgroup of a Frobenius complement must be cyclic.
Hence $U/W$ is cyclic, so
$q=4$.  In this case $C/W$ is dihedral of order $6$, so
(since $C$ contains $U$) the group $C$ is dihedral of order $12$.

Let $T'$ be the order-$3$ subgroup of $C$.  Since $T'$ is normal
in $C$, no subgroup of $C$ properly containing $T'$ can be an inertia group in
$\cC\to \cC/C$ (by Corollary~\ref{cor-ordinary}).  Thus, every orbit of $C/T'$ on
the set $\Gamma$ of fixed points of $T'$ is regular, so $|\Gamma|$ is divisible by
$4$.  Since $T$ fixes precisely two points of $\cC$, it follows that $T'$ and
$T$ are not conjugate in $\cG$, so a Sylow $3$-subgroup of $\cG$ is noncyclic,
and thus contains an elementary abelian subgroup of order $9$
\cite[23.9]{As}.
%
By Lemma~\ref{lemma-Cab}, $\cC$ has genus $q(q-1)/2=6$.  But Riemann-Hurwitz
shows that (in characteristic not $3$) an elementary abelian group of order $9$
cannot act on a genus-$6$ curve, contradiction.
\end{proof}

\begin{thm}
\label{thm-B or G}
If $\cG\ne B$ then $\cG=\SL(q)$ and $\cC\to \cC/\cG$ has precisely two branch points,
with inertia groups $B$ and $W$.
\end{thm}

\begin{proof}
Assume that $\cG \ne B$.
Consider the cover $\cC \to \cC/\cG$.  By Lemmas \ref{lemma-B-curves} and
\ref{lemma-p2}, the
inertia groups of $P_1$ and $P_2$ in this cover are $B$ and $W$,
respectively.  Since these groups are nonconjugate, $P_1$ and $P_2$
lie over distinct branch points $Q_1$ and $Q_2$.
By Lemma~\ref{lemma-Cab} and Corollary~\ref{cor-ordinary}, every
branch point of $\cC\to \cC/\cG$ has trivial second ramification group.

For a point $Q$ of $\cC/\cG$, let $\ind(Q)$ denote the sum of the different
exponents (in the cover $\cC\to \cC/\cG$) of the points lying over $Q$.
Note that $\ind(Q_1)/|\cG| = 1 - 2/|B| + |U|/|B| = 1 + (q-2)/|B|$
and $\ind(Q_2) = |\cG|$.
The Riemann-Hurwitz formula gives
\begin{align*}
q(q-1) - 2 &= -2|\cG| + \ind(Q_1) + \ind(Q_2) + \sum_{Q\notin\{Q_1,Q_2\}} \ind(Q) \\
&= (q-2)|\cG\col B| + \sum_Q \ind(Q);
\end{align*}
since any branch point $Q$ satisfies
$\ind(Q)\ge 2|\cG|/3 >q(q-1)$, it follows that $Q_1$ and $Q_2$ are the only branch
points in $\cC\to \cC/\cG$, and we must have $|\cG\col B|=q+1$.

By Lemma~\ref{MBequiv}, $T$ has index $2$ in $H:=N_{\cG}(T)$.
Thus $H$ preserves the set $\{P_1,P_3\}$ of fixed points of $T$.
Lemma~\ref{MBequiv} implies $P_1\in \cG P_3$, so
$|\cG P_3|=|\cG P_1|=|\cG\col B|=q+1$.  Since $|BP_3|=q$, it follows that
$\cG P_3=BP_3\cup\{P_1\}$.  Pick an involution $\nu\in H$.
If $\nu$ fixes $P_1$ then Lemma~\ref{lemma-B-curves}
implies $\nu\in B$; but $\nu$ must also fix $P_3$, which is
impossible since the stabilizer of $P_3$ in $B$ is $T$ (and $T$
contains no involutions).  Thus $\nu$ must swap $P_1$ and $P_3$.
By Lemma~\ref{lemma-B-curves}, $U$ is a Sylow $2$-subgroup of $\cG$;
since all involutions of $U$ are conjugate in $B$, it follows that $\nu$
is conjugate in $\cG$ to the nonidentity element of $W$, and thus fixes a
unique point of $\cG P_1$.

The orbits of $B$ on $\Lambda:=\cG P_1$ are the fixed point $P_1$
and the $q$-element orbit $BP_3$.  Since $B$ has a unique conjugacy
class of index-$q$ subgroups, this determines $\Lambda$ as a $B$-set.
The same orbit sizes occur in the action of $B$ on $\Pone(\Fq)$ induced
by the usual action  of $\PSL(q)$ on $\Pone(\Fq)$.  Thus, $\Lambda$ and
$\Pone(\Fq)$ are isomorphic $B$-sets.  We will show below that, up to
$T$-conjugacy,
there is a unique involution in the symmetric group of $\Gamma$
which normalizes $T$ and has a unique fixed point.  Since $\SL(q)$ 
contains such an involution, we can extend our isomorphism of 
$B$-sets $\Lambda\cong_B\Pone(\Fq)$ to an isomorphism of 
$\langle B,\nu\rangle$-sets, and in particular $\SL(q)$ has a subgroup
isomorphic to $\langle B,\nu\rangle$.  Since $B$ is a maximal
subgroup of $\SL(q)$, we have $\langle B,\nu\rangle\cong\SL(q)$, whence
(since $|\cG| \le |\SL(q)|$) we conclude $\cG\cong\SL(q)$.

It remains to show that, up to $T$-conjugacy, there is a unique
involution $\hat\nu$ in the symmetric group of $\Lambda$ which normalizes $T$
and has a unique fixed point.  Note that 
$T$ fixes $P_1$ and $P_3$, and $T$ is transitive on the other 
$q-1$ points of $\Lambda$.  Thus $\hat\nu$ permutes $\{P_1,P_3\}$, and the fixed point
hypothesis implies $\hat\nu$ interchanges $P_1$ and $P_3$.
Hence $\hat\nu$ fixes a unique point of $TP_3$, so we may
may identify this orbit with $T$ and assume the fixed
point is $1 \in T$. The only
order-$2$ automorphism of $T$ with no nontrivial
fixed points
is the automorphism inverting all elements of $T$,
whence $\hat\nu$ is unique up to $T$-conjugacy.
%
\end{proof}


\section{$G$-curves and hyperelliptic quotients}
\label{sec-G-curves}

Let $\bl$ be an algebraically closed field of characteristic $2$,
let $\alpha,\beta\in\bl^*$, and let $\cC:=\cC_{\alpha,\beta}$ be as in~(\ref{eq-Cab}).
We use the embedding $B\to \Aut \cC$ from Lemma~\ref{lemma-BCab}.
By Theorem~\ref{thm-B or G}, the automorphism group of
$\cC$ is either $B$ or $G:=\SL(q)$.  In this section
we determine when the latter occurs.

\begin{prop}
\label{prop-G-curves}
$\cC$ has automorphism group $G$ if and only if $\beta^2=\alpha + \alpha^2$.
\end{prop}

Set $t:=w^{q-1}$ and $y:=v/w$.  Since $T$ fixes $t$ and $y$, we
have $\bl(t,y)\subseteq \bl(v,w)^T=\bl(\cC/T)$.
Clearly $w$ has degree at most $q-1$ over $\bl(t,y)$,
and also $\bl(v,w)=\bl(y,w)$.  Thus, $\bl(\cC/T)=\bl(t,y)$.

The curve $\cC/T$ is defined by the equation
\begin{equation*}
y^q + \frac{y}{t} = \frac{\alpha + \beta}{t} + \T\Bigl(\frac{\beta}{t+1}\Bigr),
\end{equation*}
which is irreducible because $[\bl(y,t)\col \bl(t)]=q$.
Putting
\[ z:=y^2 + y + \frac{\beta}{t+1},\]  we compute
\[ \T(z)=y^q+y+\T\Bigl(\frac{\beta}{t+1}\Bigr)=
y\Bigl(1+\frac{1}{t}\Bigr)+\frac{\alpha + \beta}{t},\] and thus
\[ y=\frac{t\T(z)+ \alpha + \beta}{t+1}.\]  It follows that

\begin{lemma}
\label{lem-kummer}
$\bl(\cC/T)=\bl(t,z)$ and $\bl(\cC)=\bl(w,z)$.
\end{lemma}

Our next result gives further information about $\cC/T$.

\begin{lemma}  
$\cC/T$ is hyperelliptic of genus $q/2$, and the hyperelliptic
involution $\nu$ fixes $z$ and maps $t \mapsto (\alpha^2+\alpha+\beta^2+z)/(z^q t)$. 
\end{lemma}

\begin{proof}
Substituting our expression for $y$ (in terms of $t$ and $\T(z)$)
into the definition of $z$ gives
\begin{equation*}
z = \frac{(t\T(z) + \alpha + \beta)^2 + (t\T(z) + \alpha + \beta)(t+1) + \beta(t+1)}{(t+1)^2},
\end{equation*}
so $0 = t^2 z^q + t(\T(z)+\alpha) + (z + \alpha^2+\alpha+\beta^2)$.
By considering the order of the pole at the point
$z=\infty$ in this equation, we see that $t\notin \bl(z)$.
Thus, $[\bl(t,z)\col\bl(z)]=2$.  Our hypothesis on the ramification
in $\cC\to \cC/B$ implies that $\bl(t,z)$ has genus $q/2$.
Hence $\cC/T$ is hyperelliptic,
and the hyperelliptic involution $\nu$ fixes $z$ and maps
$t \mapsto (\alpha^2+\alpha+\beta^2+z)/(z^q t)$.
\end{proof}

Suppose in this paragraph that $\Aut_{\bbl}(\cC)\cong G$, and choose the
isomorphism so that it extends our previous embedding $B\hookrightarrow
\Aut_{\bbl}(\cC)$.  By Theorem~\ref{thm-B or G}, there are points
$P_1, P_2$ on $\cC$ whose
stabilizers in $G$ are $B$ and $W$, respectively, and moreover the
corresponding points $Q_1,Q_2$ on $\cC/G$ are the only two branch points
of $\cC\to \cC/G$.
By Lemma~\ref{MBequiv}, $H:=N_G(T)$ has order $2(q-1)$, so Lemma~\ref{lemma-gt}
implies $H$ is dihedral, hence contains $q-1$ involutions.  But all
involutions in $G$ are conjugate, and each fixes $q/2$ points of $GP_2$,
so $\cC/T\to \cC/H$ is ramified over $q/2$ points lying over $Q_2$.
Likewise, $\cC/T\to \cC/H$ is ramified over a unique point lying over $Q_1$,
so $\cC/T\to \cC/H$ has $1+q/2$ branch points and thus (since $\cC/T$ has genus
$q/2$) we find that $\cC/H$ has genus zero.
By uniqueness of the hyperelliptic involution, we must have
$\bl(\cC)^H=\bl(z)$, and each element $\mu\in H\setminus T$ is an 
involution whose restriction to $\cC/T$ is the hyperelliptic involution
$\nu$.  Now, $(w\mu(w))^{q-1}=t\rho(t)=(\alpha^2+\alpha+\beta^2+z)/z^q$ is in $\bl(z)$,
so $\bl(w\mu(w),z)/\bl(z)$ is cyclic of order dividing $q-1$; but
the dihedral group of order $2(q-1)$ has no proper normal subgroups of
even order, so $w\mu(w)\in\bl(z)$.  Thus $(\alpha^2+\alpha+\beta^2+z)/z^q$ is
a $(q-1)$-th power in $\bl(z)$, so $\beta^2=\alpha+\alpha^2$.

Conversely, we now assume that 
$\beta^2 = \alpha^2 + \alpha$ (with $\alpha\notin\F_2$, since $\beta\ne 0$).
By Lemma~\ref{lem-kummer}, there are precisely $q-1$ extensions of
$\nu$ to an embedding of $\bl(\cC)$ into its algebraic closure,
one for each $(q-1)$-th root of $\rho(t)$ (this root will be $\rho(w)$).
Since $t\rho(t) = 1/z^{q-1}$, each of these extensions
maps $w\mapsto\zeta/(zw)$ with $\zeta \in \Fq^*$ and so in particular 
leaves $\bl(\cC)=\bl(w,z)$ invariant (and thus is an automorphism
of $\bl(\cC)$).  Since $\Aut_{\bbl}(\cC)$ properly contains $B$, Theorem
\ref{thm-B or G} implies that $\Aut_{\bbl}(\cC)\cong
\SL(q)$.  This completes the proof of Proposition~\ref{prop-G-curves}.


\section{Forms of $\cC_{\alpha,\beta}$}
\label{sec-Cab}

In this section we study isomorphisms between curves of the
shape $\cC_{\alpha,\beta}$, and isomorphisms between these curves and
other curves.

\begin{prop}
\label{prop-Cab}
Let $\bl$ be an algebraically closed field of characteristic $2$.
For $\alpha,\beta,\alpha',\beta'\in\bl^*$, the curves $\cC_{\alpha,\beta}$ and 
$\cC_{\alpha',\beta'}$ are isomorphic if and only if $\alpha=\alpha'$ and $\beta=\beta'$.
\end{prop}

\begin{proof}
Let $\cC=\cC_{\alpha,\beta}$ and $\cC'= \cC_{\alpha',\beta'}$, and let $\cG=\Aut \cC$ and 
$\cG'=\Aut \cC'$.  Write the equations of $\cC$ and $\cC'$ as
$v^q+v=(\alpha+\beta)w+w^q\T(\beta/(1+w^{q-1}))$ and $(v')^q+v'=(\alpha'+\beta')w'+
(w')^q\T(\beta'/(1+(w')^{q-1}))$, respectively.
Suppose there is an isomorphism $\rho\colon \cC\to \cC'$.
Conjugation by $\rho$ induces an isomorphism 
$\theta\colon \cG\to \cG'$.  By replacing $\rho$ by its 
compositions with automorphisms of $\cC$ and $\cC'$, we can replace
$\theta$ by its compositions with arbitrary inner automorphisms of
$\cG$ and $\cG'$.

We use the embeddings $B\to \cG$ and $B\to \cG'$ from
Lemma~\ref{lemma-BCab}.  By Lemma~\ref{lemma-B-curves}, $U$ is a
Sylow $2$-subgroup of $\cG$ and $\cG'$, so (by composing $\rho$
with automorphisms) we may assume $\theta(U)=U$.
Since all index-$2$ subgroups of $U$ are conjugate under
$B$, we may assume in addition that
$\theta(H)=H$ where $H$ is a prescribed index-$2$ subgroup of $U$.  
Then $\rho$ induces an isomorphism between $\cC/U$ and $\cC'/U$ which
maps the set of branch points of $\cC/H\to \cC/U$ to the corresponding
set in $\cC'/U$.  For definiteness, choose $H$ to be the subgroup
defined in Lemma~\ref{lemma-Cab}, and choose the coordinates $w$ and
$w'$ on $\cC/U$ and $\cC'/U$.  The branch points of
each of $\cC/H\to \cC/U$ and $\cC'/H\to \cC'/U$ (in the coordinates $w$ and
$w'$) are $\{\delta : \T(\delta)=1\}\cup\{\infty\}$.

Since $B$ is the normalizer of $U$ in both $\cG$ and $\cG'$
(by Lemma~\ref{lemma-B-curves}),
it follows from $\theta(U)=U$ that
$\theta(B)=B$.  The only points of $\cC/U$ which ramify in $\cC/U\to \cC/B$
are $w=0$ and $w=\infty$, so $\rho$ must map these to $w'=0$ and 
$w'=\infty$ in some order.
Thus, $\rho(w)$ is a constant times either $w'$ or $1/w'$.  Since also
$\rho$ preserves $\{\delta : \T(\delta)=1\}\cup\{\infty\}$, we 
must have $\rho(w)=w'$.
Since $\theta(H)=H$ and the right hand side
of (\ref{eq-v}) has only simple poles, by applying $\rho$ to this equation
we see that $\alpha=\alpha'$ and $\beta=\beta'$.
\end{proof}

\begin{prop}
\label{prop-unique}
Let $k$ be a perfect field of characteristic $2$, and let $\kb$ be
an algebraic closure of $k$.  Let $\cC=\cC_{\alpha,\beta}$
where $\alpha,\beta\in \kb^*$.
Let $\cC'$ be a curve over $k$ which is isomorphic to $\cC$ over $\kb$.
Let $\ell$ be an extension of $k$ such that\/
$\Aut_\ell(\ell(\cC'))\cong\Aut_{\overline{k}}(\kb(\cC'))$.
Then:
\begin{enumerate}
\item $k$ contains\/ $\F_2(\alpha,\beta)$;
\item $\cC$ is defined over $k$; and
\item $\cC$ is isomorphic to $\cC'$ over $\ell$.
\end{enumerate}
\end{prop}

\begin{proof}
Note that $\kb(\cC) = \kb(v,w)$ where $v,w$ satisfy
\begin{equation*}
v^q + v = (\alpha+\beta)w + w^q \T\Bigl(\frac{\beta}{1+w^{q-1}}\Bigr).
\end{equation*}
If $\rho$ is any $k$-automorphism of $\kb(\cC)$, then
$\kb(\cC) = \kb(v_1,w_1)$ where $v_1:=\rho(v)$ and
$w_1:=\rho(w)$ satisfy
\begin{equation*}
v_1^q + v_1 = (\rho(\alpha)+\rho(\beta))w_1 + w_1^q \T\Bigl(\frac{\rho(\beta)}
 {1+w_1^{q-1}}\Bigr).
\end{equation*}
Thus, $\cC_{\alpha,\beta} \cong \cC_{\rho(\alpha), \rho(\beta)}$, whence (by the previous
result) $\rho$ fixes $\alpha$ and $\beta$.  Hence $\F_2(\alpha,\beta)$ is fixed by
the full group of $k$-automorphisms of $\kb(\cC)$.

By hypothesis, there is a $\kb$-isomorphism $\theta$ between
$\kb(\cC)$ and $\kb(\cC')$.  Conjugation by $\theta$ induces an isomorphism
$\Aut_k(\kb(\cC))\cong\Aut_k(\kb(\cC'))$, so in particular both of
these groups fix the same subfield of $\kb$.  Since $k$ is perfect
and $\cC'$ is defined over $k$, the subfield of $\kb$ fixed by 
$\Aut_k(\kb(\cC'))$ is just $k$, so $\F_2(\alpha,\beta)\subseteq k$.

Clearly $\cC$ is defined over $\F_2(\alpha,\beta)$, hence over $k$.
Finally, by Theorem~\ref{thm-B-curves} and 
Proposition~\ref{prop-Cab} there is an $\ell$-isomorphism
$\ell(\cC) \cong \ell(\cC')$.
\end{proof}


\section{Existence and uniqueness of polynomials}
\label{sec-polynomials}

Let $k$ be a perfect field of characteristic $2$, and let $q=2^e>2$.  
In this section we prove a preliminary version of 
Theorem~\ref{thm-intro}, in which we describe the Galois closure
of $k(x)/k(f(x))$ rather than describing the polynomials $f$.
Here $x$ is transcendental over $k$, and
we say $b,c\in k[X]$ are $k$-\emph{equivalent} if there are linear
polynomials $\ell_1,\ell_2\in k[X]$ such that $b=\ell_1\circ c\circ \ell_2$.

\begin{thm}
\label{thm-existence}
If $f\in k[X]$ is a separable polynomial of degree $(q^2-q)/2$ 
such that
\begin{enumerate}
\item the geometric monodromy group of $f$ is\/ $\SL(q)$; and
\item the extension $k(x)/k(f(x))$ is wildly ramified over at least two
places of $k(f(x))$,
\end{enumerate}
then there is a unique pair $(\alpha,\beta)\in k^*\times k^*$ with $\beta^2=\alpha+\alpha^2$ 
for which the Galois closure of $k(x)/k(f(x))$ is isomorphic to 
$(k.\Fq)(\cC_{\alpha,\beta})$.  Conversely, each such pair $(\alpha,\beta)$ actually occurs
for some $f$ with these properties, and two such polynomials are 
$k$-equivalent if and only if they correspond to the same pair $(\alpha,\beta)$.
Finally, every such $f$ is indecomposable, and
$f$ is exceptional if and only if $e$ is odd and $k \cap \Fq = \F_2$.
\end{thm}

Our proof uses a corollary of the following simple lemma (cf.\
\cite[Thm.\ 4.2A]{DM}):

\begin{lemma}
Let $G$ be a transitive permutation group on a set $\Delta$, and let $G_1$ be the
stabilizer of a point $\pi\in \Delta$.  Let $C$ be the centralizer of $G$ in
the symmetric group on $\Delta$.  Then $C\cong N_G(G_1)/G_1$, and $C$ acts faithfully
and regularly on the set of fixed points of $G_1$.  In particular, $C$ is trivial
if $G_1$ is self-normalizing in $G$.
\end{lemma}

\begin{proof}
Note that an element $\tau\in C$ is determined by the value $\tau(\pi)$ (since
$\tau(\nu(\pi))=\nu(\tau(\pi))$ for every $\nu\in G$).

If $G$ acts regularly on $\Delta$, then we can identify the action of $G$ on $\Delta$
with the action of $G$ on itself by left multiplication.  Clearly right
multiplication commutes with this action, so the map $\tau\mapsto \tau(1)$ induces
an isomorphism $C\cong G$, and $C$ acts regularly on $\Delta$.

Let $\Lambda$ be the set of fixed points of $G_1$.
Then $N_G(G_1)/G_1$ acts regularly on $\Lambda$.
Letting $\hat C$ be the centralizer of $N_G(G_1)$ in $\Sym(\Lambda)$, the previous
paragraph shows
that $\hat C\cong N_G(G_1)/G_1$ acts regularly on $\Lambda$.  Since $C$ acts on
$\Lambda$ and $C$
centralizes $N_G(G_1)$, restriction to $\Lambda$ induces a homomorphism $\theta\colon C\to \hat C$.
We see that $\theta$
is injective, since $\tau\in C$ is determined by $\tau(\pi)$.  It remains only to prove
that $\theta$ is surjective.  For $\mu\in \hat C$, $\nu\in G$ and $\lambda\in G_1$, note that
$\nu(\lambda(\mu(\pi)))=\nu(\mu(\lambda(\pi)))=\nu(\mu(\pi))$; hence the image of $\mu(\pi)$
is constant on each coset in $G/G_1$,
so the map $\nu(\pi) \mapsto \nu(\mu(\pi))$ defines a permutation $\phi$ of $\Delta$.
Plainly $\phi$
centralizes $G$ and $\theta(\phi)=\mu$, so the proof is complete.
\end{proof}

\begin{cor}
\label{cag}
Let $f\in k[X]$ be a separable polynomial, let $E$ be the Galois
closure of $k(x)/k(f(x))$, and let $\ell$ be the algebraic closure of $k$
in $E$.  Put $A:=\Gal(E/k(f(x)))$, $G:=\Gal(E/\ell(f(x)))$,
and $G_1:=\Gal(E/\ell(x))$.  If $N_G(G_1)=G_1$, then $C_A(G)=1$.
\end{cor}


\begin{proof}[Proof of Theorem~\ref{thm-existence}]
Suppose $f\in k[X]$ is a separable polynomial of degree $(q^2-q)/2$
which satisfies conditions (i) and (ii) of Theorem~\ref{thm-existence}.
Let $E$ be the Galois closure of $k(x)/k(f(x))$, and let $\ell$ be
the algebraic closure of $k$ in $E$.
Then there is an $\ell$-isomorphism between $E$ and $\ell(\cC_{\alpha,\beta})$
for some $\alpha,\beta\in\ell^*$, and also $\ell\supseteq\Fq$ (by
Corollary~\ref{cor-B} and Theorem~\ref{thm-B-curves}).
This uniquely determines the pair $(\alpha,\beta)$ (Proposition~\ref{prop-Cab}).
By Theorem \ref{thm-B or G}, the geometric monodromy group
$G:=\Gal(E/\ell(f(x)))$ equals $\Aut_\ell\ell(\cC_{\alpha,\beta})$,
so Proposition \ref{prop-G-curves} implies $\beta^2=\alpha^2+\alpha$.
By Lemma \ref{lemma-GZ} and Proposition \ref{prop-unique}, both $\alpha$ and $\beta$
are in $k$.
By Lemma \ref{lemma-gt}, the hypotheses
of the above corollary are satisfied, so no nontrivial element of $\Gal(E/k(f(x)))$
centralizes $G$.  Since every $\ell$-automorphism of $\ell(\cC_{\alpha,\beta})$
is defined over $k.\Fq$, we see that $G$ commutes with $\Gal(E/(k.\Fq)(\cC_{\alpha,\beta}))$,
so $L=(k.\F_q)(\cC_{\alpha,\beta})$.
We have proven the first sentence of Theorem \ref{thm-existence}.

Conversely, suppose $\alpha,\beta\in k^*$ satisfy $\beta^2=\alpha+\alpha^2$, and put
$\ell:=k\Fq$.
Let $E=\ell(\cC_{\alpha,\beta})$.  We have shown that $G:=\Aut_\ell E$ satisfies
$G\cong\SL(q)$, and that there
are degree-one places $P_1$ and $P_2$ of $E$ whose stabilizers in $G$ are $B$ and $W$,
respectively.  Moreover, $E$ has genus $q(q-1)/2$, and the second
ramification groups at $P_1$ and $P_2$ are trivial.  By Riemann-Hurwitz,
the only places of $E^G$ which ramify in $E/E^G$ are the places $Q_1$ and $Q_2$ which
lie under $P_1$ and $P_2$.
Let $G_1$ be a subgroup of $G$ of index $q(q-1)/2$.  Then $G_1$ is dihedral of order
$2(q+1)$, and hence contains $q+1$ involutions.  Each of the $q+1$ conjugates of $U$
contains precisely one of these involutions.  Hence there is a unique place of $E^{G_1}$
lying over $Q_1$, and its ramification index in $E/E^{G_1}$ is $2$.  Also there are precisely
$q/2$ places of $E^{G_1}$ which lie over $Q_2$ and ramify in $E/E^{G_1}$, and each has ramification
index $2$.  Thus $E^{G_1}$ has genus zero, and $Q_1$ is totally ramified in $E^{G_1}/E^G$.
Next, $A:=\Aut_k E$ satisfies $A=G.\Gal(E/k(\cC_{\alpha,\beta}))$.  Since $G$ is normal in $A$, and $G_1$ is
conjugate (in $G$) to all $(q^2-q)/2$ subgroups of $G$ having order $2q+2$, it follows that
$|N_A(G_1)\col G_1|=|\ell\col k|$ and $N_A(G_1)G=A$.  Thus, $E^{N_A(G_1)}$ is a genus-zero function
field over $k$ which contains a degree-one place that is totally ramified over
$E^A$.  We can write $E^{N_A(G_1)}=k(x)$ and $E^A=k(u)$, and by making linear fractional
changes in $x$ and $u$ we may assume that the unique place of $k(x)$ lying over
the infinite place of $k(u)$ is the infinite place.  In other words, $u=f(x)$ for some
$f\in k[X]$.  Separability of $f$ follows from
separability of $k(x)/k(u)$.  The degree of $f$ is $(q^2-q)/2$, and its geometric
monodromy group is $\SL(q)$ (since $G_1$ contains no nontrivial normal subgroup of
$\SL(q)$).  The extension $k(x)/k(f(x))$ is totally ramified over infinity, and also
is wildly ramified over another place of $k(f(x))$.

Next we show that the Galois closure of $k(x)/k(f(x))$ is $E$, or equivalently that
$N_A(G_1)$ contains no nontrivial normal subgroup of $A$.  Let $J$ be a proper normal
subgroup of $A$.  Since $G$ is normal in $A$ (and simple), $J$ must intersect $G$
trivially.  Thus each element of $J$ has shape $\nu\sigma$, where $\nu\in G$ and $\sigma\in
\Gal(E/k(\cC_{\alpha,\beta}))$ satisfy
$|\langle \nu\sigma\rangle|=|\langle\sigma\rangle|$.  In particular, $J$ is cyclic;
let $\nu\sigma$ be a generator of $J$.  Since $G$ and $J$ normalize one another and
intersect trivially, they must commute.  Write $E=\ell(v,w)$ where
$v^q+v=(\alpha+\beta)w+w^q\T(\beta/(1+w^{q-1}))$.  Let $\tau\in G$ map $(v,w)\mapsto (v+1,w)$.
Since $\tau$ commutes with both $J$ and $\sigma$, it also must commute with $\nu$.
Hence $\nu$ maps $(v,w)\mapsto (v+\alpha,w)$ for some $\alpha\in\Fq$.
For $\zeta\in\Fq^*$, let $\lambda_\zeta\in G$ map
$(v,w)\mapsto (\zeta v,\zeta w)$.  Then $\lambda_\zeta \nu \sigma(w)=\zeta w$, but
$\nu\sigma\lambda_\zeta(w)=\sigma(\zeta)w$, so $\sigma$ fixes $\zeta$.  Hence $\sigma$
fixes both $\Fq$ and $k(\cC_{\alpha,\beta})$, so it fixes $E$, whence $J=1$.  Thus the
arithmetic monodromy group of $f$ is $A$.  Since $G$ has a unique conjugacy class
of subgroups of index $(q^2-q)/2$, all of which are self-normalizing, any two
index-$(q^2-q)/2$ subgroups of $A$ which surject onto $A/G$ are conjugate.
Since $A=\Aut_k E$, it follows that there is a unique $k$-equivalence class of
polynomials $f$ which satisfy all our hypotheses for a given pair $(\alpha,\beta)$.
Conversely, $k$-equivalent polynomials have isomorphic Galois closures, hence
correspond to the same pair $(\alpha,\beta)$.  Finally, the indecomposability and
exceptionality criteria follow from Lemma~\ref{lemma-GZ}.
\end{proof}


\begin{cor}
\label{cor-existence}
There exists a separable
polynomial $f \in k[X]$ of degree $q(q-1)/2$ with two
wild branch points and 
geometric monodromy group\/ $\SL(q)$ if and only if
$k$ properly contains\/ $\F_2$.
\end{cor}

\begin{cor}
There exists a separable exceptional
polynomial $f \in k[X]$ of degree $q(q-1)/2$ with two
wild branch points and 
geometric monodromy group\/ $\SL(q)$ if and only if $e$ is odd,
$k \cap \Fq = \F_2$, and $k$ properly contains\/ $\F_2$.
\end{cor}


\section{Another nonexistence proof over $\F_2$}
\label{sec-nonexistence}

One consequence of Corollary~\ref{cor-existence} is that
there is no separable polynomial $f$
over $\F_2$ of degree $q(q-1)/2$ such that the cover
$f\colon\Pone\to\Pone$ has at least two wildly ramified branch points
and has geometric monodromy group $\SL(q)$.  In this section we
give a more direct proof of this fact, by showing that the Galois
closure of such a cover $f\colon\Pone\to\Pone$ would be a curve
having more rational points than is permitted by the Weil bound.  

\begin{thm}  
\label{thm-weil}
There is no separable polynomial $f \in \F_2[X]$
of degree $q(q-1)/2$ satisfying the following conditions:
\begin{enumerate}  
\item the geometric monodromy group of $f$ is $G:=\SL(q)$;
\item the extension\/ $\F_2(x)/\F_2(f(x))$ has precisely two 
 branch points, and in the Galois closure $E/\F_2(f(x))$
 their ramification indices are $q(q-1)$ and $2$; and 
\item all second ramification groups in $E/\F_2(f(x))$ are trivial.
\end{enumerate}
\end{thm}

\begin{remark}
By Lemma~\ref{lemma-GZ}, conditions (ii) and (iii) follow from (i)
if we assume that $f$ has two wild branch points.  Thus, the
combination of Theorem~\ref{thm-weil} and Lemma~\ref{lemma-GZ}
implies the `only if' implication in Corollary~\ref{cor-existence}.
\end{remark}

\begin{proof} 
Suppose there is an $f$ satisfying the above conditions.
The Riemann-Hurwitz formula implies that 
the genus of $E$ is $q(q-1)/2$.  

Since the two branch points of $E/\F_2(f(x))$ have nonconjugate inertia
groups, these points must be $\F_2$-rational.  Let $Q$ be the point 
with ramification index $2$.

Let $A:=\Gal(E/\F_2(f(x)))$ be the arithmetic monodromy group of $f$.  
By Corollary~\ref{cag}, $G\le A \le \Aut(G)=\SL(q).e$.
Thus, $A=G.e'$ for some $e'\mid e$.  It follows that the algebraic
closure of $\F_2$ in $E$ is $\ell:=\F_{2^{e'}}$.
Let $P$ be a place of $E$ lying over $Q$.  Let $H$ be the decomposition
group of $P$ in the extension $E/\F_2(f(x))$.
We know that the inertia group $W$ of $P$ has order $2$,
so $U:=N_G(W)$ has order $q$.
Thus, $H \le N_A(W)=\langle U, \nu\rangle$, where $\nu\in A$ has 
order $e'$ and maps to a generator of $A/G$.
Since $Q$ is $\F_2$-rational, $H/W$ surjects onto $A/G$, or 
equivalently $A=GH$.  Since $H/W$ is cyclic, it follows that 
$|H/W|$ is either $e'$ or $2e'$.

Suppose that $e'<e$.
Let $\hat\ell$ be the quadratic extension of $\ell$.  Then
$|\hat\ell| \le q$. Let $\hat P$ be a place of $\hat\ell E$ lying over $P$ (there are one or
two such places).  Since $|H/W|$ divides $[\hat\ell\col\F_2]$, the place $\hat P$
is rational
over $\hat\ell$.  Moreover, the ramification index of $\hat P$ in $\hat\ell E/\hat\ell(f(x))$\
is $2$.  Thus, $Q$ lies under $|G|/2$ rational places of $\hat\ell E$.
Since $\hat\ell E$ has genus $q(q-1)/2$, this violates the Weil bound for the
number of rational points on a curve over a finite field.

Now suppose that $e'=e$.  As noted above, $H \le N_A(W)=\langle U, \nu \rangle$.
For any $\mu \in U$, the element $(\mu \nu)^e\in H$ lies in $U$ and centralizes
$\mu \nu$, hence it centralizes $\nu$.  However, the centralizer of $\nu$
in $U$ is $W$.  Since $|\cC_U(\nu)|=2$, it follows that
no element of $N_A(W)/W$ has order $2e$, so
$H/W$ is cyclic of order $e$.
Now, as in the previous case, we obtain a contradiction by counting
points.
\end{proof}

\section{Construction of polynomials}
\label{sec-explicit}

In this section we use the results proved so far in order to
compute explicit forms of the polynomials whose existence was
proved in Theorem \ref{thm-existence}.

Let $k$ be a perfect field of characteristic $2$, and let $q=2^e>2$.
Let $\alpha,\beta\in k^*$ satisfy $\beta^2=\alpha+\alpha^2$.  

\begin{thm}
\label{thm-construct1}
The polynomial
\eq{
\hat f(X) := (\T(X)+\alpha+1)
        \prod_{\substack{\zeta^{q-1}=1\\ \zeta\ne 1}}
         \Bigl(\sum_{i=0}^{e-1}\frac{\zeta^{2^i}+\zeta}{\zeta^{2^i}+1}
          X^{2^i} + \zeta \alpha + 1\Bigr)\label{eq-f}}
is in the $k$-equivalence class corresponding to $(\alpha,\beta)$
in Theorem~\ref{thm-existence}.
\end{thm}

\begin{proof}
Let $\ell=k.\Fq$ and
$E=\ell(\cC_{\alpha,\beta})$.
Write $E=\ell(v,w)$, where $v^q+v=(\alpha+\beta)w+w^q\T(\beta/(1+w^{q-1}))$.

Let $\hat w=1/w$ and $\hat v=v^2/w+v+\beta w/(1+w^{q-1})$.
Then
\begin{align*}
\T(\hat v \hat w)&=\Bigl(\frac{v}w\Bigr)^q+\frac{v}w+\T\Bigl(\frac{\beta}{1+w^{q-1}}\Bigr) \\
 &=v\Bigl(\frac{1}w+\frac{1}{w^q}\Bigr)+\frac{\alpha+\beta}{w^{q-1}},
\end{align*}
so $k(\hat v,\hat w)=k(v,w)$.  Next,
\begin{align*}
{\hat v}^q \hat w &=
\frac{v^{2q}}{w^{q+1}}+\frac{v^q}w+\frac{\beta^qw^{q-1}}{1+w^{q^2-q}} \\
 &=\frac{v^2}{w^{q+1}}+\frac{\alpha^2+\beta^2}{w^{q-1}}+w^{q-1}\T\Bigl(\frac{\beta}{1+
w^{q-1}}\Bigr)^2
   + \frac{v}w+\alpha+\beta\\
&\qquad+w^{q-1}\T\Bigl(\frac{\beta}{1+w^{q-1}}\Bigr)
   + \frac{\beta^qw^{q-1}}{1+w^{q^2-q}} \\
 &=\frac{v^2}{w^{q+1}}+\frac{\alpha}{w^{q-1}}+\frac{w^{q-1}\beta}{1+w^{q-1}} +
\frac{v}w +\alpha+\beta \\
 &=\T(\hat v \hat w)+{\hat w}^q \hat v + \alpha.
\end{align*}
Since $[\ell(\hat v,\hat w)\col\ell(\hat w)]=[\ell(v,w)\col\ell(w)]=q$, the polynomial
$\hat w X^q+\T(\hat w X)+{\hat w}^q X + \alpha$ is irreducible over $\ell(\hat w)$, so
also $\hat v X^q+\T(\hat v X)+{\hat v}^q X+\alpha$ is irreducible over $\ell(\hat v)$.
Let $\hl=\F_{q^2}.\ell$, and $\hL=\hl.E$.  Pick $\gamma\in\hl^*$ of
multiplicative order $q+1$, and let $\delta=\gamma+1/\gamma\in\Fq^*\subseteq\ell^*$.
Let $y=(\hat v\gamma+\hat w/\gamma+1)/\delta$ and $z=(\hat v/\gamma+\hat w\gamma+1)/\delta$.
Then $\hL=\hl(\hat v,\hat w)=\hl(y,z)$.

\begin{lemma}
We have $[\hl(y,z)\col\hl(z)]=q+1$ and
\eq{\label{eq-EF}
y^{q+1}+z^{q+1}=\T(yz)+\alpha+1.}
For $\eta\in\F_{q^2}$ with $\eta^{q+1}=1$,
there is a unique element $\hat\nu_\eta\in\Aut_{\hl} \hL$ which maps
$(y,z)\mapsto (y\eta,z/\eta)$.  Moreover, $\nu_\eta:=\hat\nu_\eta|_E$
is in $\Aut_\ell E$.
\end{lemma}

\begin{proof}
We compute
\begin{align*}
y^{q+1}+z^{q+1}&=\frac{({\hat v}\gamma+\frac{{\hat w}}{\gamma}+1)(\frac{{\hat v}^q}{\gamma}+{\hat w}^q\gamma+1)}{\delta^{q+1}}
+\frac{(\frac{{\hat v}}{\gamma}+{\hat w}\gamma+1)({\hat v}^q\gamma+\frac{{\hat w}^q}{\gamma}+1)}{\delta^{q+1}}\\
&={\hat w}^q{\hat v}+{\hat w}{\hat v}^q+\frac{{\hat w}^q+{\hat v}^q+{\hat w}+{\hat v}}{\delta} \\
&=\T({\hat w}{\hat v})+\alpha+\frac{{\hat w}^q+{\hat v}^q}{\delta^q}+\frac{{\hat w}+{\hat v}}{\delta} \\
&= \T\left({\hat w}{\hat v}+\frac{{\hat w}+{\hat v}}{\delta}+\frac{{\hat w}^2+{\hat v}^2+1}{\delta^2}\right)
+\alpha+\T\Bigl(\frac{1}{\delta^2}\Bigr) \\
&= \T(yz)+\alpha+\T\Bigl(\frac{1}{\delta^2}\Bigr).
\end{align*}
Since $1/\delta=\gamma/(\gamma^2+1)=\gamma/(\gamma+1)+\gamma^2/(\gamma^2+1)$,
we have
\begin{align*}
\T\Bigl(\frac{1}{\delta^2}\Bigr)=\T\Bigl(\frac{1}{\delta}\Bigr)&=
\frac{\gamma^q}{\gamma^q+1}+\frac{\gamma}{\gamma+1} \\
&=\frac{\frac{1}{\gamma}}{\frac{1}{\gamma}+1}+\frac{\gamma}{\gamma+1}=1.
\end{align*}
Since $\hL=\hl(y,z)$ has genus $q(q-1)/2$ where $y$ and $z$ satisfy
equation (\ref{eq-EF}) of total degree $q+1$, this equation
must define a smooth
(projective) plane curve, and in particular must be irreducible.  Thus
$[\hl(y,z)\col\hl(z)]=q+1$.  Now existence and uniqueness of $\hat\nu_\eta$
are clear.  A straightforward computation yields that $\hat\nu_\eta$ maps
\begin{align*}
{\hat w}\mapsto& \frac{1}{\delta^2}\left(\delta+\Bigl(\frac{\gamma}{\eta}+
\frac{\eta}{\gamma}\Bigr)+{\hat w}\Bigl(\frac{\eta}{\gamma^2}+\frac{\gamma^2}{\eta}\Bigr)
+{\hat v}\Bigl(\eta+\frac{1}{\eta}\Bigr)\right) \\
{\hat v}\mapsto& \frac{1}{\delta^2}\left(\delta+\Bigl(\gamma\eta+\frac{1}{\gamma
\eta}\Bigr)+{\hat w}\Bigl(\eta+\frac{1}{\eta}\Bigr)+{\hat v}\Bigl(\gamma^2\eta+
\frac{1}{\gamma^2\eta}\Bigl)\right).
\end{align*}
Since $\ell(\hat\nu_\eta({\hat v}),\hat\nu_\eta({\hat w}))=E$, it follows that
$\hat\nu_\eta$ induces an automorphism of~$E$.
\end{proof}

We now compute the subfield of $E$ fixed by an index-$(q^2-q)/2$ subgroup of 
$G:=\Aut_\ell E\cong\SL(q)$.
There is a unique element
$\tau\in\Aut_\ell E$ such that $\tau\colon (\hat v,\hat w)\mapsto (\hat w,\hat v)$.
Note that $\tau$ maps $(y,z)$ to $(z,y)$, and the group
$G_1:=\langle \tau,\{\nu_\eta : \eta^{q+1}=1\}\rangle$ is dihedral
of order $2q+2$.
Hence the subfield of $\hL$
fixed by $G_1$ contains $\hl(yz)$.  Multiplying equation (\ref{eq-EF}) by
$y^{q+1}$, we see that $[\hl(y,z)\col\hl(yz)]=[\hl(y,yz)\col\hl(yz)]\le 2q+2$,
so $\hl(yz)$ is the subfield of $\hL$ fixed by $G_1$.
Moreover, since
\begin{equation*}
yz=\left(\frac{\hat v+\hat w}\delta\right)^2+\frac{\hat v+\hat w}\delta+\hat v\hat w+\frac{1}{\delta^2}
\end{equation*}
lies in $E$, the subfield of $E$ fixed by $G_1$ is $\ell(yz)$.

Next we compute an invariant of $G$.
Recall that $\SL(q)$ can be written as $CTU$, where $T$ is the diagonal
subgroup, $U$ is a unipotent subgroup, and $C$ is a cyclic subgroup of order $(q+1)$.
We can choose $U$ to be the set of maps $\sigma_\xi\colon (v,w)\mapsto(v+\xi,w)$ with $\xi\in\Fq$,
so $E^U=\ell(w)$.
We can choose $T$ to be the set of maps
$\mu_\zeta\colon (v,w)\mapsto(\zeta^{-1}v,\zeta^{-1}w)$ with $\zeta^{q-1}=1$, and
$C$ to be the set of maps $\nu_\eta$ defined in the above lemma.
Hence the product
\begin{equation*}
\prod_{\eta^{q+1}=1}\prod_{\zeta^{q-1}=1}\prod_{\xi\in\Fq}
 \nu_\eta\mu_\zeta\sigma_\xi\Bigl(\frac{\delta}{w}+1\Bigr)\end{equation*}
is $G$-invariant.  Since this product is the $q$-th power of
\begin{equation*}
u:=\prod_{\eta^{q+1}=1}\prod_{\zeta^{q-1}=1} \nu_\eta\mu_\zeta\Bigl(
\frac{\delta}{w}+1\Bigr),
\end{equation*}
also $G$ fixes $u$.  Since $1/w=(z\gamma+1+y\gamma^{-1})/\delta$, we have
\begin{equation*}
u = 
  \prod_{\zeta^{q-1}=1}\prod_{\eta^{q+1}=1}
       (\eta\zeta y\gamma^{-1} + \zeta + 1 + \eta^{-1}\zeta\gamma z).
\end{equation*}

By the following lemma,
\begin{align*}
u&=\prod_{\zeta^{q-1}=1}\prod_{\eta^{q+1}=1}
       (\eta\zeta y + \zeta + 1 + \eta^{-1}\zeta z) \\
&= (y^{q+1}+z^{q+1})\prod
          _{\substack{\zeta^{q-1}=1\\ \zeta\ne 1}}
       \biggl(\zeta^2(y^{q+1}+z^{q+1}) \,\,+ \\
&
\qquad\qquad\qquad\qquad\qquad\qquad (\zeta^2+1)
       \Bigl(1+\T\Bigl(yz\frac{\zeta^2}{\zeta^2+1}\Bigr)\Bigr)\biggr) \\
&= (\T(yz)+\alpha+1)\prod
          _{\zeta\in\Fq\smallsetminus\F_2}
       \biggl(\zeta\bigl(\T(yz)+\alpha+1\bigr) \,\,+\\
&\qquad\qquad\qquad\qquad\qquad\qquad\qquad (\zeta+1)
       \Bigl(1+\T\Bigl(yz\frac{\zeta}{\zeta+1}\Bigr)\Bigr)\biggr) \\
&= (\T(yz)+\alpha+1)\prod
          _{\zeta\in\Fq\smallsetminus\F_2}
         \Bigl(\sum_{i=0}^{e-1}\frac{\zeta^{2^i}+\zeta}{\zeta^{2^i}+1}(yz)^{2^i}
          + \zeta \alpha + 1\Bigr).
\end{align*}
Thus $u=\hat f(yz)$ where $\hat f$ is the polynomial defined in~(\ref{eq-f}).
It follows that $[\ell(yz)\col\ell(u)]=\deg(\hat f)=(q^2-q)/2$.
Since $E^{G_1}=\ell(yz)$ and
$E^G\supseteq\ell(u)$ and $[E^{G_1}\col E^G]=(q^2-q)/2$,
it follows that $E^G=\ell(u)$.  Now, $G_1$ contains no nontrivial
normal subgroup of $G$, so $E$ is the Galois closure of
$\ell(yz)/\ell(u)$, whence $G$ is the geometric monodromy group of $\hat f$.
Clearly $\hat f$ is fixed by $\Gal(\ell/k)$, so $\hat f\in k[X]$.
By Theorem~\ref{thm-B or G}, the extension $E/E^G$ has two wildly
ramified branch points, so Theorem \ref{thm-existence} implies that
$\hat f$ is in the $k$-equivalence class corresponding to the pair $(\alpha,\beta)$.
\end{proof}

\begin{lemma} 
The following identity holds in $k[Y,Z]$:
\begin{equation*}
\prod_{\omega^{q+1}=1} (\omega Y + 1 + \omega^{-1}Z) = Y^{q+1} + Z^{q+1} + \T(YZ) + 1.
\end{equation*}
\end{lemma}

\begin{proof}
By applying the transformation $(Y,Z)\mapsto(\omega Y,Z/\omega)$,
we see that $\prod (\omega Y+1+\omega^{-1}Z)-Y^{q+1}-Z^{q+1}$ is a polynomial
$h(YZ)\in k[YZ]$, with degree at most $q/2$ and constant term 1.
If we substitute $Y=Z=\omega/(\omega^2+1)$ (where $\omega^{q+1}=1$ and
$\omega \ne 1$), we see that $YZ=\omega^2/(\omega^4+1)$ is a root of $h$.
These roots of $h$ are precisely the trace~1 elements of $\F_q$,
namely the roots of $\T(YZ)+1$.  Hence $h(YZ)$ and $\T(YZ)+1$ have the same
roots and the same constant term, and $\T(YZ)+1$ is squarefree with
$\deg(\T+1)\ge\deg(h)$, so $h(YZ)=\T(YZ)+1$.
\end{proof}

\begin{remark}
Once one knows `where to look' for these polynomials -- especially,
what should be the Galois closure $E$ of $k(x)/k(f(x))$ -- one can
give direct proofs of their properties.  But such proofs would seem unmotivated,
since we know no way to guess what $E$ should be besides appealing
to the results in this paper.
\end{remark}


\section{Another form for the polynomials}
\label{sec-explicit2}

In the previous section we computed the polynomials
whose existence was proved in Theorem \ref{thm-existence}.
Our expression for the polynomials was concise, but involved a product.
In this section we prove Theorem~\ref{thm-intro} and
Corollary~\ref{cor-intro} by writing the polynomials without any
sums or products other than the usual $\T(X)=X^{q/2}+X^{q/4}+\dots+X$.
Here $q=2^e>2$ and $k$ is a perfect field of characteristic $2$.
Also $\alpha,\beta\in k^*$ satisfy $\beta^2=\alpha^2+\alpha$.

\begin{thm} \label{final}
The expression
\[
f(X):=\left(\frac{\T(X)+\alpha}{X}\right)^q\cdot\left(\T(X)+\frac{\T(X)+\alpha}{\alpha+1}\cdot
\T\Bigl(\frac{X(\alpha^2+\alpha)}{(\T(X)+\alpha)^2}\Bigr)\right)
\]
defines a polynomial which lies in the $k$-equivalence class corresponding to
$(\alpha+1,\beta)$ in Theorem~\ref{thm-existence}.
\end{thm}

\begin{proof}
First we show that $f$ is a polynomial.  Writing $h(X):=X^q f(X)$, we have
\begin{align*}
h &= (\T(X)+\alpha)^q \cdot\T(X) + \frac{(\T(X)+\alpha)^{q+1}}{\alpha+1}\cdot
\T\Bigl(\frac{X(\alpha^2+\alpha)}{(\T(X)+\alpha)^2}\Bigr) \\
&= (\T(X)+\alpha)^q \cdot\T(X) + \frac{1}{\alpha+1}\sum_{i=0}^{e-1}
X^{2^i}(\alpha^2+\alpha)^{2^i}(\T(X)+\alpha)^{q+1-2^{i+1}}.
\end{align*}
Thus $h$ is a polynomial divisible by $X\cdot (\T(X)+\alpha)$, and
moreover $h$ is monic of degree $q(q+1)/2$.
We now determine the multiplicity of $X$ as a divisor of $h$.
This multiplicity is unchanged if we replace $h$ by
\[
\hat h := h\cdot\Bigl(h+\frac{(\T(X)+\alpha)^{q+1}}{\alpha+1}\Bigr);
\]
writing $c:=X(\alpha^2+\alpha)/(\T(X)+\alpha)^2$,
we compute
\begin{align*}
\hat h &= h^2 + h\cdot \frac{(\T(X)+\alpha)^{q+1}}{\alpha+1} \\
&= (\T(X)+\alpha)^{2q}\cdot \T(X)^2 + \T(X)\cdot\frac{(\T(X)+\alpha)^{2q+1}}{\alpha+1} +
\frac{(\T(X)+\alpha)^{2q+2}}{\alpha^2+1}\cdot\T(c^2+c).
\end{align*}
Substituting $\T(c^2+c)=c^q+c$, and reducing mod $X^{2q}$, we find that
\begin{align*}
\hat h &\equiv \alpha^{2q} \T(X)^2 + \frac{\alpha^{2q}}{\alpha+1} (\T(X)^2+\alpha\T(X)) +
\frac{(\T(X)+\alpha)^2}{\alpha^2+1}X^q(\alpha^2+\alpha)^q \\
&\qquad+
\frac{\alpha^{2q}}{\alpha^2+1}X(\alpha^2+\alpha)
 \pmod{X^{2q}} \\
&= \frac{\alpha^{2q}}{\alpha+1}\left( \T(X)^2(\alpha+1) + \T(X)^2+\alpha\T(X) + \alpha X\right) +
\frac{(\T(X)+\alpha)^2}{\alpha^2+1}X^q(\alpha^2+\alpha)^q
 \\
&= \frac{\alpha^q}{\alpha^2+1} X^q \left(
\alpha^{q+1}(\alpha+1) + (\T(X)+\alpha)^2(\alpha+1)^q \right),
\end{align*}
so $X^q$ divides $\hat h$, whence $f$ is a polynomial divisible
by $(\T(X)+\alpha)$.
Furthermore, $X$ divides $f$ (equivalently $X^{q+1}$ divides $\hat{h}$)
precisely when $\alpha \in \F_q$, in which case
$X^2$ exactly divides $f$.  Since $h$ is monic of degree $q(q+1)/2$, it
follows that $f$ is monic of degree $q(q-1)/2$.

We now show that $f/(\T(X)+\alpha)$ is in $k[X^2]$.
It suffices to show that $\bar f:=X^q f/(\T(X)+\alpha)$ is in $k[X^2]$.  We compute
\begin{align*}
\bar f &= \T(X)(\T(X)+\alpha)^{q-1} + 
\frac{(\T(X)+\alpha)^q}{\alpha+1}\cdot
\T\Bigl(\frac{X(\alpha^2+\alpha)}{(\T(X)+\alpha)^2}\Bigr) \\
&=
(\T(X)+\alpha)^q + \alpha(\T(X)+\alpha)^{q-1} + 
\frac{1}{\alpha+1}\sum_{i=0}^{e-1} (X(\alpha^2+\alpha))^{2^i} (\T(X)+\alpha)^{q-2^{i+1}}.
\end{align*}
The summands with $i>0$ are polynomials in $X^2$.
Thus, there exists $b\in k[X]$ such that
\begin{align*}
\bar f &= b(X^2) + \alpha(\T(X)+\alpha)^{q-1} + \alpha X(\T(X)+\alpha)^{q-2} \\
&= b(X^2) + \alpha(\T(X)+\alpha)^{q-2}(\T(X)+\alpha+X),
\end{align*}
so indeed $\bar f\in k[X^2]$, whence $f/(\T(X)+\alpha)$ is in $k[X^2]$.

By Theorem~\ref{thm-construct1}, the polynomial
\[
\hat f(X):=(\T(X)+\alpha)\prod_{\substack{\zeta^{q-1}=1 \\ \zeta\ne 1}}
\Bigl(\sum_{i=1}^{e-1}\frac{\zeta^{2^i}+\zeta}{\zeta^{2^i}+1}X^{2^i}+\zeta(\alpha+1)+1\Bigr)
\]
is in the $k$-equivalence class corresponding to $(\alpha+1,\beta)$ in
Theorem~\ref{thm-existence}.
By Lemma~\ref{lemma-GZ}, the extension $k(x)/k(\hat f(x))$ has precisely
two branch points; one of these points is totally ramified, and the ramification
index at any point of $k(x)$ lying over the other branch point is at most $2$.
Since $k(x)/k(\hat f(x))$ is totally ramified over the infinite place, there is
a unique finite branch point.
But plainly $\hat f(X)=(\T(X)+\alpha)\hat b(X)^2$ for some nonconstant
$\hat b\in k[X]$, so $\hat f(x)=0$ is the finite branch point, and thus
$\hat b(X)$ is squarefree and coprime to $(\T(X)+\alpha)$.
We will show that every root $\delta$ of $\hat f$ is a root of $f$;
it follows that the multiplicity of $\delta$ as a root of $f$ is at least as big as
the corresponding multiplicity for $\hat f$.
Since $f$ and $\hat f$ have the same degree and the same leading coefficient, we
conclude that $f=\hat f$.

It remains to prove that every root of $\hat f$ is a root of $f$.
Recall that, in the function field $\kb(y,z)$ where
$y^{q+1}+z^{q+1}=\T(yz)+\alpha$, we have the identity
\begin{align*}
\hat f(yz) &= \prod_{\zeta^{q-1}=1}\prod_{\eta^{q+1}=1}(\eta\zeta y+\zeta+1+
\frac{\zeta}{\eta}z) \\
&= (y^{q+1}+z^{q+1})\prod_{\substack{\zeta^{q-1}=1\\ \zeta\ne 1}}\prod_{\eta^{q+1}=1}
(\eta\zeta y+\zeta+1 + \frac{\zeta}{\eta}z).
\end{align*}
Let $\delta$ be a root of $\hat{f}$.  Pick $\hE\in\kb^*$ and $\hF\in\kb$
such that
$\delta=\hE\hF$ and $\hE^{q+1}+\hF^{q+1}=\T(\hE\hF)+\alpha$: such $\hE,\hF$ exist
because substituting $\hF=\delta/\hE$ into the latter equation (and clearing denominators)
gives a polynomial in $\hE$ which is not a monomial, and thus has a nonzero
root.
If $\T(\delta)=\alpha$ then we already know that $f(\delta)=0$.
If $\delta=0$ then $\hF=0$ and $\hE^{q+1}=\alpha$, so
\begin{align*}
0=\hat{f}(0) &= \prod_{\zeta^{q-1}=1}\prod_{\eta^{q+1}=1}(\eta\zeta\hE+\zeta+1)\\
&= \prod_{\zeta^{q-1}=1}(\zeta^{q+1}\hE^{q+1} + (\zeta+1)^{q+1}) \\
&= \prod_{\zeta^{q-1}=1}(\zeta^2\alpha + \zeta^2+1) \\
&= (\alpha+1)^{q-1} +1.
\end{align*}
Thus $\alpha\in\F_q$, so $X^2$ divides $f$.

Henceforth we assume $\alpha\ne\T(\delta)$ and $\delta\ne0$.
This implies $\eta\zeta\hE + \zeta + 1 + \hF\zeta/\eta=0$ for some
$\zeta,\eta$ with $\zeta\in\F_q\setminus\F_2$ and $\eta^{q+1}=1$.
By replacing $\hE$ and $\hF$
with $\eta\hE$ and $\hF/\eta$, we may assume $\eta=1$, so
\[
\hF = \hE + 1+\frac{1}{\zeta}.
\]
Write $\hz:=1+1/\zeta$, and note that $\hz\in\F_q\setminus\F_2$.
Since $\delta=\hE\hF$, we compute
\begin{align*}
\T(\delta)+\alpha &= \hE^{q+1}+\hF^{q+1} \\
&= \hE^{q+1} + \hE^{q+1} + \hz\hE^q +
   \hz^q\hE + \hz^{q+1} \\
&= \hz\hE^q + \hz\hE + \hz^2
\end{align*}
and
\begin{align*}
\T(\delta) &= \T(\hE^2+\hz\hE) \\
&= \hz\hE^q + \hz\hE + \T(\hE^2+\hz^2\hE^2).
\end{align*}
Thus
\[
\alpha+\hz^2 = \T(\hE^2+\hz^2\hE^2),
\]
so
\[
\sqrt{\alpha}+\hz = \T(\hE+\hz\hE).
\]
Adding the last two equations gives
\[
\alpha+\sqrt{\alpha}+\hz^2+\hz = \hE^q + \hE + \hz\hE^q + \hz\hE = (1+\hz)(\hE^q+\hE),
\]
so
\begin{align*}
\T(\delta)+\alpha &= \hz^2 + \hz(\hE^q+\hE) \\
&= \hz^2 + \frac{\hz}{1+\hz}
 \left(\alpha+\sqrt{\alpha}+\hz^2+\hz\right)  \\
&= \frac{\hz}{1+\hz}\left(\alpha+\sqrt{\alpha}\right)
\end{align*}
and
\begin{align*}
\T\Bigl(\frac{\delta}{\hz^2}\Bigr) &= \T\Bigl(\frac{\hE^2}{\hz^2}+\frac{\hE}{\hz}\Bigr)
\\ &= \frac{\hE^q}{\hz} + \frac{\hE}{\hz} \\
&= 1 + \frac{\alpha+\sqrt{\alpha}}{\hz^2+\hz}.
\end{align*}

Writing $\tilde f(X):=X^q f(X)/(\T(X)+\alpha)^q$, we have
\begin{align*}
\tilde f(\delta) &=
\T(\delta) + \frac{\T(\delta)+\alpha}{\alpha+1}\cdot \T\Bigl(\frac{\delta(\alpha^2+\alpha)}
{(\T(\delta)+\alpha)^2}\Bigr) \\
&= \frac{\alpha+\hz\sqrt{\alpha}}{1+\hz} + \frac{\hz\sqrt{\alpha}}{(1+\hz)(\sqrt{\alpha}+1)}\cdot
\T\Bigl(\frac{\delta(\alpha^2+\alpha)(1+\hz)^2}{\hz^2(\alpha^2+\alpha)}\Bigr) \\
&= \frac{\alpha+\hz\sqrt{\alpha}}{1+\hz}
 + \frac{\hz\sqrt{\alpha}}{(1+\hz)(\sqrt{\alpha}+1)}\cdot \T\Bigl(\delta +\frac{\delta}{\hz^2}\Bigr) \\
&= \frac{\alpha+\hz\sqrt{\alpha}}{1+\hz}
 + \frac{\hz\sqrt{\alpha}}{(1+\hz)(\sqrt{\alpha}+1)}\cdot \frac{\hz+\alpha+\sqrt{\alpha}(1+\hz)}
{\hz} \\
&= \frac{\alpha+\hz\sqrt{\alpha}+\sqrt{\alpha}(\sqrt{\alpha}+\hz)}{1+\hz},
\end{align*}
so $\tilde f(\delta)=0$ and thus $f(\delta)=0$, which completes the proof.
\end{proof}

\begin{remark}
The above proof is not completely satisfying, since it is a verification
that $f(X)$ has the desired property, rather than a derivation of the
simple expression for $f(X)$.  We do not have a good explanation why
the polynomial in Theorem~\ref{thm-existence} can be written in such a
simple form.
\end{remark}

We conclude the paper by proving the results stated in the introduction.

\begin{proof}[Proof of Theorem~\ref{thm-intro}]
In case $k$ is perfect, the result follows from Theorem~\ref{thm-existence}
and Theorem~\ref{final}.  For general $k$, let $\tilde{k}$ denote the
perfect closure of $k$.  Let $f\in k[X]$ satisfy properties
(i) and (ii) of Theorem~\ref{thm-intro}.  Then $f$ satisfies
the same properties over the perfect field $\tilde{k}$, so
$f$ is $\tilde{k}$-equivalent to $f_{\alpha}$ for some
$\alpha\in\tilde{k}\setminus\F_2$.  We will show that this implies
$f$ is $k$-equivalent to $f_{\alpha}$, and that $\alpha\in k$.
Since the monodromy groups of $f$ over $k$ are the same as those
over $\tilde{k}$, indecomposability and exceptionality
of $f$ over $k$ are equivalent to the corresponding properties over
$\tilde{k}$.  Since $\tilde{k}\cap\F_q = k\cap\F_q$ (because
$\F_q/\F_2$ is separable), the result follows.

It remains to prove that if $f(X):=\delta+\eta f_{\alpha}(\zeta X+\gamma)$
is in $k[X]$, where $\delta,\eta,\alpha,\zeta,\gamma\in\tilde{k}$
with $\eta\zeta\ne 0$ and $\alpha\notin\F_2$, then $\delta,\eta,\alpha,
\zeta,\gamma$ are in $k$.  The terms of $f_{\alpha}(X)$ of degree
at least $(q^2-3q)/2$ are $\T(X) X^{(q^2-2q)/2}+\alpha X^{(q^2-3q+2)/2}$
and (if $q=4$) $(\alpha+1)X^2$.  Hence the coefficients of
$X^{q^2/2-q+2}$ and $X^{q^2/2-q+1}$ in $f(X)$ are $\eta \zeta^{q^2/2-q+2}$ and
$\eta\zeta^{q^2/2-q+1}$, and since these are in $k^*$, we must have
$\zeta,\eta\in k^*$.
The coefficients of $X^{(q^2-3q+2)/2}$ and $X^{(q^2-2q)/2}$ in $f(X)$ are
$\alpha \eta \zeta^{(q^2-3q+2)/2}$ and $\eta \zeta^{(q^2-2q)/2} \T(\gamma)$,
so $\alpha\in k^*$ and $\T(\gamma)\in k$, whence
$\T(\gamma)^2+\T(\gamma)=\gamma^q+\gamma$ is in $k$.
The coefficient of $X^{(q^2-3q)/2}$ in $f(X)$ is
$\eta \zeta^{(q^2-3q)/2}(\alpha \gamma+\gamma^q)$
(plus $\eta(\alpha+1)\zeta^2$ if $q=4$), so $\gamma$ is in $k$.
Finally, we conclude that $\delta=f(0)-\eta f_{\alpha}(\gamma)$ is in $k$.
\end{proof}

\begin{proof}[Proof of Corollary~\ref{cor-intro}]
First assume $f\in k[X]$ is a separable indecomposable exceptional polynomial
in case (i) of Theorem~\ref{fgsthm}.  Then the geometric monodromy group
$G$ of $f$ is solvable, and the degree $d$ of $f$ is prime and not equal to
$p$.  By \cite[Thm.~4]{Mu-Schur}, it follows that $f$ is $\bar{k}$-equivalent
to either $X^d$ or $D_d(X,1)$.  By \cite[Lemma~1.9]{T}, $f$ is $k$-equivalent
to either $X^d$ or $D_d(X,a)$ with $a\in k^*$.  These polynomials $f(X)$ are
separable and indecomposable.  We verify exceptionality by examining the
factorization of $f(X)-f(Y)$ in $\bar{k}[X,Y]$, given for instance in
\cite[Prop.~1.7]{T}.

Now consider case (iii) of Theorem~\ref{fgsthm}.
In this case, Corollary~\ref{cor-intro} for $p=3$ is \cite[Thm.~1.3]{GZ}.
So suppose $p=2$, and
let $f\in k[X]$ be a separable indecomposable exceptional polynomial
of degree $d=q(q-1)/2$ where $q=2^e>2$ with $e>1$ odd.
By Theorem~\ref{fgsthm}, the arithmetic monodromy group $A$ of $f$
is $\PGammaL(q)$, and thus $G$ has a transitive normal subgroup isomorphic
to $\PSL(q)$.
The desired result follows from \cite[Thm.~4.3]{GZ} if
$\bar{k}(x)/\bar{k}(f(x))$ has no finite branch points, or if the
Galois closure $E$ of this extension does not have genus $(q^2-q)/2$.
If neither of these conditions hold, then \cite[Thm.~2.1]{GZ}
implies that $G=\PSL(q)$ and $E/\bar{k}(f(x))$ has precisely one
finite branch point, whose inertia group has order $2$ and whose
second ramification group is trivial.  In particular, $f$ satisfies
conditions (i) and (ii) of Theorem~\ref{thm-intro}, so in this
case the result follows from Theorem~\ref{thm-intro}.
\end{proof}


\end{document}